%
%
\input amstex 
\documentstyle{conm-p}
\NoBlackBoxes

\issueinfo{00}
  {}
  {}
  {2003}

\define\bbR{{\bold R}}
\define\bbC{{\bold C}}

\define\rto{\text{\bf R}\hskip-.5pt^2}

\define\ve{\hskip.7pt\varepsilon}
\define\ah{\alpha}
\define\gm{\gamma}
\define\da{\delta}
\define\de{\delta}
\define\om{\omega}
\define\ly{{\hskip.4ptl}}
\define\py{p}
\define\qe{q}

\define\xp{E}
\define\r{,\hskip2pt}

\define\lap{\Lambda\hskip-1pt^+\hskip-1.4ptM}
\define\lam{\Lambda\hskip-1pt^-\hskip-1.4ptM}
\define\lapm{\Lambda\hskip-1pt^\pm\hskip-1.4ptM}
\def\tm{{T\hskip-.3ptM}}

\def\tmc{{[\tm]^\bbC}}

\def\biv{{[\tm]^{\wedge2}}}
\def\bivx{{[T_xM]^{\wedge2}}}
\define\e{\Cal E}
\define\tv{\Cal V}
\define\x{\Cal X}
\define\y{\Cal Y}
\define\z{\Cal Z}

\define\ef{P}
\define\fe{F}
\define\zz{z\hskip-.5ptz}

\define\w{W}
\define\wy{\w^{\text{\fivebf(\hskip-.7pt+\hskip-.7pt)}}} 
\define\diml{-di\-men\-sion\-al}
\define\fdi{fi\-nite-di\-men\-sion\-al}
\define\ifd{in\-fin\-ite-di\-men\-sion\-al}
\define\dimr{\dim_{\hskip.4pt\bbR\hskip-1.2pt}}
\define\dimc{\dim_{\hskip.4pt\bbC\hskip-1.2pt}}
\define\spanc{\text{\rm Span}_{\hskip.5pt\bbC\hskip.3pt}}
\define\inv{-in\-var\-i\-ant}

\define\sky{skew-sym\-me\-try}
\define\dia{di\-ag\-o\-nal\-izable}
\define\diy{di\-ag\-o\-nal\-iza\-bil\-i\-ty}
\define\rc{real/com\-plex}
\define\cli{com\-plex-lin\-e\-ar} 
\define\cbi{com\-plex-bi\-lin\-e\-ar} 
\define\cdi{com\-plex-di\-ag\-o\-nal\-izable}
\define\cdy{com\-plex-di\-ag\-o\-nal\-iza\-bil\-i\-ty}
\define\vs{vector space}
\define\rvs{real vector space}
\define\cvf{complex vector field}
\define\vf{vector field}
\define\kf{Killing field}
\define\la{Lie algebra}
\define\lig{Lie group}
\define\os{open subset}

\define\vb{vector bundle}
\define\tb{tangent bundle}
\define\rctb{real/com\-plex\-i\-fied tangent bundle}
\define\mf{manifold}
\define\mfd{-man\-i\-fold}
\define\psr{pseu\-do\hs-Riem\-ann\-i\-an}
\define\psrm{pseu\-do-Riem\-ann\-i\-an manifold}
\define\lcc{Levi-Civita connection}
\define\ch{cur\-va\-ture-ho\-mo\-ge\-ne\-ous}
\define\lh{locally homogeneous}
\define\nbd{neighborhood}
\define\ci{$\,C^\infty$}
\define\of{$\,1$-form}
\define\mppp{\hbox{$-$\hskip1pt$+$\hskip1pt$+$\hskip1pt$+$}}
\define\mmpp{\hbox{$-$\hskip1pt$-$\hskip1pt$+$\hskip1pt$+$}}
\define\mmmp{\hbox{$-$\hskip1pt$-$\hskip1pt$-$\hskip1pt$+$}}
\define\pppp{\hbox{$+$\hskip1pt$+$\hskip1pt$+$\hskip1pt$+$}}
\define\mpmp{\hbox{$-$\hskip1pt$\pm$\hskip1pt$+$}}
\define\mpmpp{\hbox{$-$\hskip1pt$\pm$\hskip1pt$+$\hskip1pt$+$}}
\define\mmpmp{\hbox{$-$\hskip1pt$-$\hskip1pt$\pm$\hskip1pt$+$}}

\define\hs{\hskip.7pt}
\define\nh{\hskip-.7pt}
\define\ptmi{\phantom{i}}
\define\ptmii{\phantom{ii}}
 

\define\a{}
\define\f{\thetag}
\define\ff{\tag}

\def\bes{1}
\def\bkv{2}
\def\bra{3}
\def\buo{4}
\def\but{5}
\def\bva{6}
\def\cle{7}
\def\cpa{8}
\def\caw{9}
\def\hdg{10\hs}
\def\fkm{11}
\def\jen{12\hs}
\def\kmc{13}
\def\kpr{14\hs}
\def\pet{15\hs}
\def\sit{16}
\def\tak{17}
\define\id{0}
\define\pr{1}
\define\op{2}
\define\ec{3}
\define\co{4}
\define\cl{5}
\define\cn{6}
\define\ri{7}
\define\fb{8}
\define\ut{9}
\define\cf{10}
\define\ms{11}
\define\li{12}
\define\rf{13}
\define\tr{A}
\define\sn{1}
\define\cu{2}
\define\cv{3}
\define\bg{4}
\define\nr{5}
\define\rt{6}
\define\lw{7}
\define\vi{8}
\define\nx{9}
\define\nw{10}
\define\ea{11}
\define\vw{12}
\define\av{13}
\define\xj{14}
\define\lm{15}
\define\ls{16}
\define\pa{17}
\define\th{18}
\define\fy{19}
\define\dw{20}
\define\fw{21}
\define\fn{22}
\define\nd{23}
\define\lz{24}
\topmatter
\title Cur\-va\-ture\hs-ho\-mo\-ge\-ne\-ous indefinite 
Einstein metrics in dimension four\hs: 
the di\-ag\-o\-nal\-izable case\endtitle
\author Andrzej Derdzinski\endauthor
\leftheadtext{ANDRZEJ DERDZINSKI}%
\rightheadtext{CUR\-VA\-TURE\hs-HO\-MO\-GE\-NE\-OUS INDEFINITE EINSTEIN 
METRICS}%

\address Department of Mathematics, Ohio State University, Columbus, Ohio 43210
\endaddress


\email andrzej\@math.ohio-state.edu\endemail





\endtopmatter

\document
\voffset=0pt
\hoffset=54.5pt
\document 
\head\S\id. Introduction\endhead 
A \psrm\ $\,(M,g)\,$ is called {\it \ch\/} if the algebraic type of its 
metric/curvature pair $\,(g,R)\,$ is the same at all points, i.e., if for 
any $\,x,y\in M\hs$ some isomorphism $\,T_xM\to T_yM\hs$ sends $\,g(x),R(x)\,$ 
to $\,g(y),R(y)$. 

Every \lh\ \psr\ manifold is, obviously, \ch. The converse proposition fails; 
counterexamples with positive-definite metrics were first found by Takagi 
\cite{\tak} and, on compact manifolds, by Ferus, Karcher and M\"unzner 
\cite{\fkm}; see also \cite{\bkv}. Analogous examples with indefinite metrics 
have been known even longer (\cite{\bra} -- \cite{\cle}, \cite{\kmc}).

The present paper provides a classification, up to local isometries, of all 
those \ch\ \psr\ four\mfd s $\,(M,g)\,$ which are Einstein and have, at some 
(or every) point $\,x$, a {\it \cdi} curvature operator 
$\,R(x):\bivx\to\bivx$. (The last condition means that the \cli\ 
extension of $\,R(x)\,$ to the complexification of the bivector space 
$\,\bivx$ is \dia.) It turns out that all such \mf s are \lh\ and, 
in fact, either locally symmetric, or locally isometric to a \lig\ with a 
left\inv\ indefinite metric of a specific type; see Theorems \a\cl.1, \a\cn.1 
and \a\ri.1. In those theorems we assume constancy of eigenvalues of the 
curvature operator, which sounds weaker than 
cur\-va\-ture-ho\-mo\-ge\-ne\-i\-ty, but, in the \cdi\ case, is actually 
equivalent to it; cf.\ \cite{\hdg}, p.~701 and Remark~6.19 on p.~472.

The metric $\,g\,$ can have any signature. Using a sign change, we may assume 
that $\,g\,$ is {\it Riem\-ann\-i\-an}\hs, {\it neutral\/} or {\it 
Lo\-rentz\-i\-an}\hs, that is, has one of the sign patterns
$$\pppp\,\,,\qquad\mmpp\,\,,\qquad\mppp\,\,.\ff\sn$$
Two known families of \ch\ Einstein four\mfd s, one Lo\-rentz\-i\-an 
\cite{\bra} and one neutral (\cite{\hdg}, p.~705), give rise to \ifd\ 
spaces of local-isometry types. By contrast, for the \mf s classified here, 
the analogous space is clearly \fdi\ (see above). Also, our \diy\ assumption 
always holds for {\it Riemannian\/} manifolds, as the curvature operator is 
self-adjoint, and for Riemannian metrics our theorem becomes the result of 
\cite{\hdg}, mentioned below. In the Lo\-rentz\-i\-an case, the \cdy\ 
condition means that the curvature is of the Petrov type I at each point, cf.\ 
\cite{\hdg}, p.~659.

Some types of \ch\ Einstein four\mfd s have already been classified. This 
includes locally symmetric spaces (\cite{\cpa}, \cite{\caw}; cf.\ \cite{\hdg}, 
pp.~662--663); Brans's classification \cite{\bra} of Lo\-rentz\-i\-an Einstein 
metrics representing the Pe\-trov type III at every point (a condition that 
implies cur\-va\-ture-ho\-mo\-ge\-ne\-i\-ty); as well as the Riemannian case 
(\cite{\hdg}, Corollary~7.2, p.~476), in which the metrics in question are all 
locally symmetric (see also \S\ri).

The text is organized as follows. In sections \op\ -- \co\ we introduce our 
``model spaces'', using a construction basically due to Petrov \cite{\pet}. 
The classification result is stated in sections \cl\ -- \ri\ and then proven 
in sections \fb\ -- \rf.

\head\S\pr. Preliminaries\endhead
Our conventions about the curvature tensor $\,R=R^\nabla$ of any connection 
$\,\nabla\,$ in a \rc\ \vb\ $\,\e\hs$ over a \mf\ $\,M$, its Ricci tensor 
$\,\hs\text{\rm Ric}\hs\,$ when $\,\e\hs$ is the tangent bundle $\,\tm$, and 
the scalar curvature $\,\hs\text{\rm s}\hs\,$ in the case where $\,\nabla\,$ 
is the \lcc\ of a given \psr\ metric $\,g\,$ on $\,M$, are
$$\alignedat2
&\text{\rm\ptmi i)}\quad&&
R(u,v)\psi\,=\,\nabla_{\!v}\nabla_{\!u}\psi\,-\,\nabla_{\!u}\nabla_{\!v}\psi\,
+\,\nabla_{[u,v]}\psi\,,\\
&\text{\rm ii)}\quad&&
\text{\rm Ric}\,(u,w)\,=\,\,\text{\rm Trace}\,[v\mapsto R(u,v)w]\,,\qquad
\text{\rm s}\,\,=\,\,\text{\rm Trace}_g\hs\text{\rm Ric}\,,\endalignedat
\ff\cu$$
for any (local) $\,C^2$ sections $\,u,v,w\,$ of $\,\tm\,$ and $\,\psi\,$ of 
$\,\e$.

A \psr\ manifold $\,(M,g)\,$ with $\,\dim M=n\,$ is said to be an {\it 
Einstein manifold\/} \cite{\bes} if $\,n\ge3\,$ and 
$\,\hs\text{\rm Ric}\,=\,\text{\rm s}\hskip1.2ptg/n$, while, if $\,n\ge4$, 
formulae 
$\,\sigma=\,\text{\rm Ric}\,-\hs(2n-2)^{-1}\,\text{\rm s}\hskip1.2ptg\,$
and $\,\w\,=\,\hs R\hs\,-\,(n-2)^{-1}\hskip1.2ptg\wedge\sigma\,$ define the 
{\it Schouten tensor\/} $\,\sigma\,$ and {\it Weyl tensor\/} $\,\w\hs$ of 
$\,(M,g)$. Here $\,\wedge\,$ is the exterior product of \of s valued in \of s, 
obtained using the valuewise multiplication which is also provided by 
$\,\wedge\hs$, so that the result is a $\,2$-form valued in $\,2$-forms.

For $\,(M,g)\,$ as above, we denote  $\,\biv$ the vector bundle of bivectors 
over $\,M$, with the fibres $\,\bivx$, $\,x\in M$. There exists a unique \psr\ 
fibre metric $\,\langle\,,\rangle\,$ in $\,\biv\,$ such that 
$\,\langle v\wedge u\hs,v'\wedge u'\rangle=g(v,v')\hskip1.2ptg(u,u')
-g(v,u')\hskip1.2ptg(u,v')\,$ for any $\,x\in M\,$ and $\,v,u,v',u'\in T_xM$.
Both $\,R\hs,\w\hs$ are four-times covariant tensor fields on $\,M\,$ sharing 
the (skew)sym\-me\-try properties of the curvature tensor, which allows us to 
treat them as morphisms acting on bivectors and self-adjoint relative to 
$\,\langle\,,\rangle\,$ at each point; in this way, $\,R\,$ gives rise to the 
{\it curvature operator}
$$R:\biv\to\biv\hskip12pt\text{\rm with}\hskip10pt
\langle R(u\wedge v),w\wedge w'\hs\rangle\,=\,\hs g(R(u,v)w,w'\hs)\ff\cv$$
for $\,x\in M\,$ and $\,u,v,w,w'\in T_xM$. When $\,(M,g)\,$ is four\diml\ and 
oriented, another important morphism $\,\biv\to\biv$ is the {\it Hodge star\/} 
$\,*\hs$, given by 
$\,*(e_1\wedge e_2)\hs=\,\ve_3\ve_4\hskip1.4pte_3\wedge e_4$ for any 
$\,x\in M\,$ and any positive-oriented orthonormal basis $\,e_1,\dots,e_4$ of 
$\,T_xM$, where $\,\ve_a=g(e_a,e_a)\in\{1,-\hs1\}\,$ (no summation). This 
well-known description of $\,*\,$ (cf. \cite{\hdg}, formula~(37.13) on p.~639) 
is clearly equivalent to its more common definition 
$\,\ah\wedge\beta\,=\,\langle*\ah\hs,\hs\beta\,\rangle\,\text{\rm vol}\hs\,$ 
for any bivectors $\,\ah,\beta$, where $\,\hs\text{\rm vol}\hs\,$ is the {\it 
volume four-vector}, equal to $\,e_1\wedge\ldots\wedge\hskip1pte_4$ for any 
$\,e_1,\dots,e_4$ as above.

Let $\,(M,g)\,$ be an oriented \psr\ $\,4$\mfd. Then $\,[\hs\w,\hs*\hs]=0\hs$, 
that is, the morphisms $\,\w\hs,\hs*:\biv\to\biv$ commute (cf.\ \cite{\sit}), 
while our formula for $\,*\,$ gives $\,*\hs^2=\,\text{\rm Id}\,$ 
if $\,g\,$ is Riemannian ($\pppp$) or neutral ($\mmpp$), and 
$\,*\hs^2=-\,\text{\rm Id}\hs\,$ when $\,g\,$ is Lo\-rentz\-i\-an ($\mppp$). 
In the Lo\-rentz\-i\-an case, this turns $\,\biv$ into a {\it complex vector 
bundle\/} of fibre dimension $\,3$, in which $\,*\,$ is the multiplication by 
$\,\,i\hs$, and, as $\,[\hs\w,\hs*\hs]=0\hs$, the Weyl tensor $\,\w\hs$ is a 
{\it \cli\ bundle morphism\/} $\,\biv\to\biv$. In the Riemannian and 
neutral cases, the self-adjoint involution $\,*:\biv\to\biv$ gives rise to the 
orthogonal decomposition $\,\biv=\lap\oplus\lam$, where $\,\lapm$, the 
$\,(\pm\hs1)$-eigenspace bundles of $\,*\hs$, are real vector bundles of fibre 
dimension $\,3$, called the bundles of {\it self-dual\/} and {\it 
anti-self-dual\/} bivectors in $\,(M,g)$. As $\,[\hs\w,\hs*\hs]=0\hs$, both 
$\,\lapm$ are $\,\w$\inv, which leads to the restrictions 
$\,\w^\pm:\lapm\to\lapm$ of $\,\w$, called the {\it self-dual\/} and {\it 
anti-self-dual\/} Weyl tensors of $\,(M,g)$. See \cite{\sit} and \cite{\hdg}, 
pp.~637 -- 651.
\remark{Remark \a\pr.1}For a \psr\ {\it Einstein\/} manifold $\,(M,g)\,$ of 
dimension $\,n\ge4$, the difference $\,R-\w:\biv\to\biv$ of the morphisms 
$\,R\hs,\w\,$ clearly equals the constant $\,\hs\text{\rm s}/[n(n-1)]\,$ times 
$\,\hs\text{\rm Id}\,=(g\wedge g)/2$. If $\,n=4\,$ and $\,M\,$ is oriented, 
relation $\,[\hs\w,\hs*\hs]=0\,$ thus gives $\,[R\hs,\hs*\hs]=0\hs$, i.e., in 
the Lo\-rentz\-i\-an case the curvature operator $\,R:\biv\to\biv$ is \cli, 
while in the Riemannian and neutral cases both $\,\lapm$ are $\,R$\inv; we 
will call the restriction $\,R^{\hs+}:\lap\to\lap\,$ of $\,R\,$ the {\it 
self-dual curvature operator\/} of $\,(M,g)$.
\endremark
\remark{Remark \a\pr.2}Let $\,x\,$ be a point in an oriented \psr\ 
$\,4$\mfd\ $\,(M,g)\,$ having one of the sign patterns \f{\sn}, and let 
$\,u\in T_xM\,$ be a vector such that $\,g(u,u)\ne0\,$ and the subspace 
$\,u\wedge u^\perp$ of $\,\bivx$ formed by all $\,u\wedge v\,$ with 
$\,v\in u^\perp$ is invariant under the Weyl tensor $\,W(x):\bivx\to\bivx$. 
Then
\widestnumber\item{(b)}\roster
\item"(a)" In the Riemannian and neutral cases, the restriction  
$\,u\wedge u^\perp\to\Lambda^{\!+\!}_xM\,$ of the orthogonal projection 
$\,\bivx\to\Lambda^{\!+\!}_xM\,$ is a linear isomorphism under which 
$\,W(x):u\wedge u^\perp\to u\wedge u^\perp$ corresponds to 
$\,\w^+(x):\Lambda^{\!+\!}_xM\to\Lambda^{\!+\!}_xM$. 
\item"(b)"In the Lo\-rentz\-i\-an case, the real subspace $\,u\wedge u^\perp$ 
spans $\,\bivx$ as a complex \vs, and the Weyl tensor 
$\,\w(x):\bivx\to\bivx$ is the unique 
\cli\ extension of $\,\w(x):u\wedge u^\perp\to u\wedge u^\perp$.
\endroster
In fact, our formula for $\,*\,$ applied to $\,e_1,\dots,e_4$ with $\,u=re_1$ 
for some $\,r>0\,$ shows that $\,\Cal H=u\wedge u^\perp$ and its 
$\,*\hs$-image $\,*\Cal H\,$ together span $\,\bivx$, and so, for dimensional 
reasons, $\,\Cal H\cap*\Cal H=\{0\}$. This gives (b). Now let $\,g\,$ be 
Riemannian or neutral. As $\,\Cal H\cap*\Cal H=\{0\}$, the space $\,\Cal H\,$ 
contains no nonzero (anti)self-dual bivectors. The projection 
$\,\bivx\to\Lambda^{\!+\!}_xM$, which has the kernel $\,\Lambda^{\!-\!}_xM$, 
is therefore injective on $\,\Cal H$, i.e., constitutes an isomorphism 
$\,\Cal H\to\Lambda^{\!+\!}_xM$. Finally, since $\,\lapm$ are $\,\w$\inv, the 
projection commutes with $\,\w(x)$, and (a) follows.
\endremark

\head\S\op. One particular family of metrics\endhead
The construction described here goes back to Petrov; see \cite{\pet}, p.~185. 

Let $\,\x\,$ be a real \vs\ of any dimension $\,n\ge3\,$ with a 
co\-di\-men\-sion-one subspace $\,V\subset\x\,$ and an element 
$\,u\in\x\smallsetminus V$, and let $\,\langle\,,\rangle\,$ be a nondegenerate 
symmetric bilinear form in $\,V$. If a linear operator $\,\fe:V\to V\,$ is 
{\it self-adjoint\/} relative to $\,\langle\,,\rangle$, that is, 
$\,\langle\fe v,v'\rangle=\langle v,\fe v'\rangle\,$ for all $\,v,v'\in V$, 
then, choosing any $\,\da\in\{1,-\hs1\}$, we define a Lie-algebra 
multiplication $\,[\hskip2.5pt,\hskip1pt]\,$ in $\,\x\,$ and a nondegenerate 
symmetric bilinear form $\,g\,$ in $\,\x\,$ by
$$\alignedat2
&\text{\rm\ptmi i)}\quad&&
[u,v]\,=\,\fe v\,,\qquad[v,v']\,=\,0\qquad\text{\rm whenever}\quad
v,v'\in V,\\
&\text{\rm ii)}\quad&&
g(u,u)\hs=\hs\da\hs,\quad g(u,v)\hs=\hs0\hs,\quad g(v,v')\hs
=\hs\langle v,v'\rangle\hskip11pt\text{\rm for\ all}\hskip8pt
v,v'\in V.\endalignedat\ff\bg$$
Let there also be given an $\,n$\diml\ real \mf\ $\,M\,$ such that $\,\x$, 
rather than being just an abstract \la, is a {\it simply transitive \la\ of 
\vf s} on $\,M$, as defined below in the appendix. An explicit description of 
such $\,M\,$ is given in the last paragraph of this section; another option is 
to choose $\,M\,$ to be the underlying \mf\ of a \lig\ $\,G\,$ whose \la\ of 
left\inv\ \vf s is isomorphic to $\,\x$. Formula \f{\bg.ii} now defines a 
\psr\ metric $\,g\,$ on $\,M\,$ such that $\,g(u,v)\,$ is constant whenever 
$\,u,v\in\x$, or, in Lie-group terms, $\,g\,$ is invariant under left 
translations in $\,G$. If $\,\nabla,\hs R\,$ and $\,\hs\text{\rm Ric}\hs\,$ 
denote the \lcc, curvature tensor and Ricci tensor of this metric $\,g$, and 
$\,v,v',w\in V\,$ are treated, along with $\,u\hs$, as \vf s on $\,M$, then
$$\alignedat2
&\text{\rm a)}\hskip8pt&&
\nabla_{\!u}u\,=\,\nabla_{\!u}v\,=\,0\,,\quad
\nabla_{\!v}u\,=\,-\,\fe v\,,\quad
\nabla_{\!v}w\,=\,\da\langle\fe v,w\rangle\hs u\,,\\
&\text{\rm b)}\hskip8pt&&
R(u,v)u\,=\,-\,\fe^2v\,,\quad R(v,w)u\,=\,0\,,
\quad R(u,w)v\,=\,\da\langle\fe^2w,v\rangle\hs u\,,\\
&\text{\rm c)}\hskip8pt&&
R(v,v')w=\da\langle\fe v',w\rangle\hs\fe v-\da\langle\fe v,w\rangle\hs\fe v',
\hskip11pt\text{\rm d)}\hskip5pt
R(u\wedge v)=-\hs\da u\wedge\fe^2v\hs,\endalignedat\ff\nr$$
with $\,R(u\wedge v)\,$ as in \f{\cv}. Namely, by \f{\bg} with 
$\,\langle\fe v,v'\rangle=\langle v,\fe v'\rangle$, the connection 
$\,\nabla\,$ in $\,\tm\,$ {\it defined\/} by \f{\nr.a} is torsionfree and 
$\,\nabla g=0\hs$, so that $\,\nabla\,$ must be the \lcc\ of $\,g$, while 
\f{\nr.b\r c} follow from \f{\cu.i}, \f{\nr.a} and \f{\bg.i}, and \f{\nr.d} is 
clear since, by \f{\cv} and \f{\nr.b\r c}, 
$\,R(u\wedge v)\hs+\hs u\wedge\fe^2v\,$ is orthogonal to $\,u\wedge w\,$ and 
$\,v'\wedge w\,$ for $\,v',w\in V\hs$ (cf.\ \S\pr), and hence to all bivectors 
at every point. Also,
$$\text{\rm Ric}\,(u,u)=-\hs\text{\rm Trace}\,\fe^2\nh,\hskip7pt
\text{\rm Ric}\,(u,v)\,=\,0\hs,\hskip6pt\text{\rm Ric}\,(v,w)
=-\hs\da\langle\fe v,w\rangle\,\text{\rm Trace}\,\fe,\ff\rt$$
for $\,v,w\in V\nh$. In fact, by \f{\cu.ii} and \f{\bg.ii}, 
$\,\,\text{\rm Ric}\,(u,u)$, 
$\,\hs\text{\rm Ric}\,(v,w)\,-\,\da g(R(u,v)u,w)\,$ and 
$\,\,\text{\rm Ric}\,(u,v)\,$ are the traces of the operators $\,V\to V\hs$ 
given by $\,v\mapsto R(u,v)u\hs$, $\,v'\mapsto R(v,v')w$, and 
$\,w\mapsto\,\text{\rm pr}\hs[R(u,w)v]$, where $\,\hs\text{\rm pr}:\x\to V\,$ 
is the orthogonal projection, so that \f{\rt} is immediate from \f{\nr.b\r c}. 

An $\,n$\diml\ \mf\ $\,M\,$ admitting a simply transitive \la\ $\,\x\,$ of 
\vf s (cf.\ the appendix) with a vector subspace $\,V\subset\x\,$ and a linear 
operator $\,\fe:V\to V\,$ such that $\,\dim V=n-1\,$ and the Lie bracket in 
$\,\x\,$ satisfies \f{\bg.i} for some $\,u\in\x\smallsetminus V$, can be 
constructed as follows. We fix $\,V\hs$ and $\,\fe:V\to V\nh$, then set 
$\,M=V\times(0,\infty)\,$ and let $\,\x=\,V+\bbR u\,$ be the space of \vf s on 
$\,M\,$ spanned by the linear \vf\ $\,u\,$ with $\,u(x,t)=(-\fe x,t)\,$ for 
$\,(x,t)\in M$, and by $\,V\nh$, each $\,v\in V\,$ being identified with the 
constant \vf\ $\,(v,0)$. Now \f{\bg.i} follows since $\,[v,w]=d_vw-d_wv\,$ 
for \vf s $\,v,w\,$ on any \os\ $\,\,U\,$ of a \fdi\ \vs\ $\,\x$, treated as 
functions $\,\,U\to\x$. 

\head\S\ec. The Einstein case\endhead
Let $\,\x,n,V,u,\fe,\langle\,,\rangle,\da,M\,$ have the properties listed in 
\S\op, and let $\,g\,$ be the \psr\ metric with \f{\bg.ii} on the $\,n$\diml\ 
\mf\ $\,M$. Then $\,g\,$ is Einstein if and only if one of the following 
conditions holds:
\widestnumber\item{(iii)}\roster
\item"(i)"$\fe\,$ equals some real scalar $\,\lambda\,$ times the identity.
\item"(ii)"$\fe\ne0\hs$, while $\,\,\text{\rm Trace}\,\fe=0\,$ and 
$\,\fe^2=0\hs$.
\item"(iii)"$\fe^2\ne0\,$ and 
$\,\,\text{\rm Trace}\,\fe=\,\text{\rm Trace}\,\fe^2=0\hs$.
\endroster
To see this, note that each of (i) -- (iii) implies, by \f{\rt}, that 
$\,g\,$ is Einstein. Conversely, let $\,g\,$ be Einstein; then either 
$\,\,\text{\rm Trace}\,\fe\ne0\,$ (and hence \f{\rt} for $\,v,w\,$ yields 
(i)), or $\,\,\text{\rm Trace}\,\fe=0\,$ and so, by \f{\rt}, 
$\,\,\text{\rm Trace}\,\fe^2=0\hs$, which in turn gives (i) (when $\,\fe=0$), 
or (ii) (when $\,\fe\ne0\,$ and $\,\fe^2=0$), or, finally, (iii) (when 
$\,\fe^2\ne0$).
\proclaim{Lemma \a\ec.1}For\/ $\,\x,n,V,u,\fe,\langle\,,\rangle,\da,M\,$ and\/ 
$\,g\,$ as above, suppose that\/ $\,g\,$ is an Einstein metric, so that we 
have\/ {\rm(i)}, {\rm(ii)} or {\rm(iii)}.

In case\/ {\rm(i)}, $\,g\,$ has the constant sectional curvature\/ 
$\,-\hs\da\lambda^2$.

In case\/ {\rm(ii)}, under the additional assumption that\/ $\,n=4$, the 
metric\/ $\,g\,$ is flat.

In case\/ {\rm(iii)}, $\,g\,$ is Ricci-flat but\/ {\rm not} locally symmetric. 
\endproclaim
In fact, the assertion about (i) follows from \f{\nr.b\r c}. Next, if 
$\,n=4$, (ii) gives $\,\fe(V)\subset\,\text{\rm Ker}\,\fe\ne V\,$ and 
$\,\dim\hs[\fe(V)]+\hs\dim\hs[\text{\rm Ker}\,\fe]=\dim V=3$, i.e., 
$\,\dim\hs[\fe(V)]=1\,$ and $\,\dim\hs[\text{\rm Ker}\,\fe]=2$. Thus, the 
right-hand sides in \f{\nr.b\r c} both vanish: the former since 
$\,\fe^2=0\hs$, the latter in view of \sky\ in $\,\fe v,\fe v'\in\fe(V)\,$ 
with $\,\dim\hs[\fe(V)]=1$. This proves our claim about (ii).

Finally, in case (iii), $\,\hs\text{\rm Ric}\,=0\,$ by \f{\rt}, 
while $\,(\nabla_{\!w}R)(u,v)v'=\nabla_{\!w}[R(u,v)v']-
R(\nabla_{\!w}u,v)v'-R(u,\nabla_{\!w}v)v'-R(u,v)\nabla_{\!w}v'$ for 
$\,v,v',w\in V\nh$, and so 
$\,\da(\nabla_{\!w}R)(u,v)v'=
-\hs\langle\fe^2v,v'\rangle\hs\fe w+\langle\fe v,v'\rangle\hs\fe^2w
-\langle\fe^2w,v'\rangle\hs\fe v+\langle\fe w,v'\rangle\hs\fe^2v\,$ by 
\f{\nr}. If $\,R\,$ were parallel, setting $\,w=v\,$ and applying 
$\,g(\,\cdot\,,v'')\,$ with any $\,v''\in V\hs$ (see \S\pr) we would get 
$\,\fe v\wedge\fe^2v=0\,$ for all $\,v\in V\nh$. Every 
$\,\fe v\in\fe(V)\smallsetminus\{0\}\,$ thus would be an eigenvector of 
$\,\fe\nh$, making $\,\fe^2$ a multiple of $\,\fe\nh$, contrary to (iii) (cf.\ 
Remark \a\ec.2).\quad\qed
\remark{Remark \a\ec.2}If $\,\hs\text{\rm Trace}\,\fe=0\,$ for an operator 
$\,\fe:V\to V\,$ in a \fdi\ \vs\ $\,V\hs$ and $\,\fe^2$ equals a nonzero 
scalar times $\,\fe$, then $\,\fe=0\hs$. In fact, let $\,r\fe^2=2\fe\,$ and 
$\,r\ne0\hs$. Then $\,A^2=\,\text{\rm Id}\hs\,$ for the operator 
$\,A=\,\hs\text{\rm Id}\hs-r\fe\,$ in $\,V$, and so 
$\,A=\pm\,\text{\rm Id}\hs\,$ on some subspaces $\,V_\pm$ with 
$\,V=V_+\oplus V_-$. Hence $\,\hs\text{\rm Trace}\,A=n_+\nh-n_-$, where 
$\,n_\pm=\dim V_\pm$, while $\,\hs\text{\rm Trace}\,A=\dim V=n_+\nh+n_-$ as 
$\,A=\,\hs\text{\rm Id}\hs-r\fe\,$ and $\,\hs\text{\rm Trace}\,\fe=0\hs$. 
Thus, $\,n_-=0\hs$, i.e., $\,V=V_+$, $\,A=\,\text{\rm Id}\hs\,$ and 
$\,\fe=0\hs$.
\endremark

\head\S\co. The curvature operator\endhead
Given a fixed sign $\,\pm\hs$, formulae
\widestnumber\item{(b)}\roster
\item"(a)"$V=\bbC\times\bbR\,$ and 
$\,\langle(z,t),(z',t')\rangle=\,\text{\rm Im}\,\zz'\pm\,tt'$ for 
$\,(z,t),(z',t')\in V$,
\item"(b)"$\fe(z,t)=(\py\qe z,\py t)$, with $\,\qe=e^{2\pi i/3}$ and 
any fixed $\,\py\in\bbR\smallsetminus\{0\}$,
\endroster
define a \rvs\ $\,V\hs$ with $\,\dim V=3$, a nondegenerate symmetric bilinear 
form $\,\langle\,,\rangle\,$ in $\,V\hs$ with the sign pattern $\,\mpmp\hs$, 
and a self-adjoint operator $\,\fe:V\to V$ satisfying (iii) in \S\ec. (See 
also the last paragraph of this section.)

In fact, $\,\langle\fe(z,t),(z',t')\rangle
=\langle(\py\qe z,\py t),(z',t')\rangle
=\,\py\hs(\text{\rm Im}\,\qe zz'\pm tt')\,$ is symmetric in $\,(z,t),(z',t')$, 
while (iii) holds for $\,\fe\,$ since $\,\fe^2(z,t)=(\py^2\qe^2z,\py^2t)$, 
$\,\hs\qe=(\sqrt{3\,}i-1)/2\,$ and $\,\qe^2=\qe^{-1}=\overline \qe$. 
\remark{Remark \a\co.1}Let $\,B:V\to V\,$ be a linear operator in an 
$\,n$\diml\ \rvs\ $\,V$. As in \S\id, we call $\,B\,$ {\it \cdi} if its \cli\ 
extension $\,B:V^{\hs\bbC}\to V^{\hs\bbC}$ to the complexification of $\,V\hs$ 
is \dia. Clearly, \ (a)\hskip7ptif $\,B\,$ is \dia, it is \cdi; \ 
(b)\hskip7pt$B\,$ is \cdi\ whenever its characteristic polynomial has $\,n\,$ 
distinct complex roots; \ (c)\hskip7ptif $\,V\hs$ is the underlying real space 
of a complex \vs\ in which $\,B\,$ acts \cli\-ly, then \cdy\ of $\,B\,$ is 
equivalent to \diy\ of $\,B\,$ as a \cli\ operator.
\endremark
\example{Example \a\co.2}Let a four\mfd\ $\,M\,$ and an indefinite metric 
$\,g\,$ on $\,M\,$ be chosen as in \S\op\ using $\,n=4$, some 
$\,\da\in\{1,-\hs1\}$, and $\,V,\langle\,,\rangle,\fe\,$ defined in (a), (b) 
above for any fixed sign $\,\pm\,$ and $\,\py\in\bbR\smallsetminus\{0\}$. 
According to Lemma \a\ec.1, $\,g\,$ is Ricci-flat but not locally symmetric. 
By \f{\bg.ii}, the sign pattern of $\,g\,$ is $\,\mpmpp\,$ (when $\,\da=1$) 
or $\,\mmpmp\,$ (when $\,\da=-\hs1$). We consider two cases:
\widestnumber\item{(ii)}\roster
\item"(i)"$\da=1\,$ and the sign $\,\pm\,$ is $\,+\hs$. Thus, $\,g\,$ is a 
Ricci-flat Lo\-rentz\-i\-an metric.
\item"(ii)"$\da=\mp\hs1$, so that $\,g\,$ is a neutral $\,(\mmpp)\,$ 
Ricci-flat metric.
\endroster
In both cases, $\,(M,g)\,$ is locally homogeneous and, by \f{\bg.ii}, locally 
isometric to a Lie group with a left\inv\ metric. (See Corollary \a\tr.3 in 
the appendix.)

Also, the curvature operator $\,R\,$ of $\,(M,g)\,$ is \cdi\ at every point. 
Namely, by \f{\nr.d}, $\,R\,$ leaves invariant the subbundle $\,\Cal H\,$ of 
$\,\biv$ spanned by all $\,u\wedge v\,$ with $\,v\in V$. Also, again by 
\f{\nr.d}, $\,R:\Cal H\to\Cal H\,$ is, at every point, algebraically 
equivalent to $\,-\hs\da\fe^2:V\to V$. On the other hand, $\,\fe^2$ is \cdi\ 
by Remark~\a\co.1(b), since $\,\fe^2/\py^2$ has the characteristic roots 
$\,1,\qe,\overline \qe$, and our assertion follows, in case (i), from 
Remark~\a\pr.2(b) for any fixed orientation of $\,M$, combined with 
Remark~\a\co.1(c), and, in case (ii), from Remark~\a\pr.2(a) applied to both 
orientations of $\,M$.
\endexample
\remark{Remark \a\co.3}As we just saw, for $\,(M,g)\,$ obtained in Example 
\a\co.2, the curvature operator $\,R\,$ (case (i)), or its self-dual 
restriction $\,R^{\hs+}\nh$, for either orientation (case (ii)), has the 
complex eigenvalues $\,\lambda,\lambda e^{2\pi i/3},\lambda e^{4\pi i/3}$ with 
$\,\lambda\in\bbR\smallsetminus\{0\}$. Also, for every locally symmetric \psr\ 
Einstein $\,4$\mfd s with a \cdi\ curvature operator, $\,R\,$ (or, 
$\,R^{\hs+}$) has a multiple eigenvalue (\cite{10}, pp.~662--663). Finally, 
according to sections \cl\ -- \ri\ below, Example \a\co.2 describes, locally, 
all possible $\,4$\diml\ \ch\ \psr\ Einstein \mf s with the sign patterns 
\f{\sn}, which are not locally symmetric. Thus, the algebraic types of 
curvature operators realized by \ch\ \psr\ Einstein $\,4$\mfd s are quite 
special, in analogy with the result of \cite{\kpr} for \ch\ {\it Riemannian\/} 
\mf s of dimension $\,4$.
\endremark
\medskip
The claim made in the three lines following (a), (b) above remains valid if 
one replaces (b) with $\,\fe(z,t)=(\pm\hs it,\,\text{\rm Re}\,z)$. This leads, 
as in Example \a\co.2, to another \lh\ Ricci-flat \psr\ $\,4$\mfd\ $\,(M,g)$, 
except that, for similar reasons, its curvature operator is {\it not\/} \cdi.

\head\S\cl. A classification theorem for the Lorentz\-i\-an case\endhead
In the following theorem, proven in \S\rf, the \diy\ assumption about the 
curvature operator amounts to its \cdy; see Remark \a\co.1(c). 
\proclaim{Theorem \a\cl.1}Let\/ $\,(M,g)\,$ be an oriented four\diml\ 
Lorentz\-i\-an Einstein manifold whose curvature operator, treated as a \cli\ 
\vb\ morphism\/ $\,R:\biv\to\biv\nh$, is \dia\ at every point and has complex 
eigenvalues that form constant functions\/ $\,M\to\bbC\hs$. Then\/ $\,(M,g)\,$ 
is \lh, and one of the following three cases occurs\/{\rm:}
\widestnumber\item{(c)}\roster
\item"(a)"$(M,g)$\hskip2.55ptis a space of constant curvature.
\item"(b)"$(M,g)$\hskip2.55ptis locally isometric to the Riemannian product of 
two \psr\ surfaces having the same constant Gaussian curvature.
\item"(c)"$(M,g)$\hskip2.55ptis locally isometric to Petrov's Ricci-flat \mf\ 
of\/ {\rm Example \a\co.2(i)}.
\endroster
Furthermore, $\,(M,g)\,$ is locally symmetric in cases\/ {\rm(a)} -- {\rm(b)}, 
but not in\/ {\rm(c)}, and in case\/ {\rm(c)} it is locally isometric to a 
\lig\ with a left\inv\ metric.
\endproclaim

\head\S\cn. A classification theorem in the neutral case\endhead
The next theorem will be proven at the end of \S\rf. For the definitions of 
$\,R^{\hs+}$ and \cdy, see Remarks \a\pr.1 and \a\co.1. 
\proclaim{Theorem \a\cn.1}Let the self-dual curvature operator\/ 
$\,R^{\hs+}:\lap\to\lap\,$ of an oriented four\diml\ Einstein \mf\/ 
$\,(M,g)\,$ of the metric signature\/ $\,\mmpp\,$ be \cdi\ at every point, 
with complex eigenvalues forming constant functions\/ $\,M\to\bbC\hs$. If\/ 
$\,\nabla\nh R^{\hs+}\ne0\,$ somewhere in\/ $\,M$, then\/ $\,(M,g)\,$ is \lh, 
namely, locally isometric to a \lig\ with a left\inv\ metric.

More precisely, $\,(M,g)\,$ then is locally isometric to one of Petrov's 
Ricci-flat \mf s, described in\/ {\rm Example \a\co.2(ii)}.
\endproclaim
Theorem \a\cn.1 sounds much stronger than its Riemannian analogue, i.e.,  
Theorem \a\ri.2 in \S\ri: an assumption about $\,R^{\hs+}$ yields a complete 
local description of the metric in the former result, but only an assertion 
about $\,R^{\hs+}$ in the latter. However, if the clause 
``$\hs\nabla\nh R^{\hs+}\ne0\,$ somewhere\hs'' were to be included among the 
hypotheses of Theorem \a\ri.2, as it is in Theorem \a\cn.1, the conclusion 
of Theorem \a\ri.2 would become an equally strong nonexistence statement.

\head\S\ri. The Riemannian case\endhead
For Riemannian metrics, our assertion amounts to the following theorem, in 
which the assumption of \cdy\ is redundant, as the curvature operator is 
self-adjoint at every point; cf.\ Remark \a\co.1(a). See also \cite{\jen}.
\proclaim{Theorem \a\ri.1\ {\rm(\cite{\hdg}, Corollary~7.2 on p.~476)}}If the 
curvature operator of a four\diml\ Riemannian Einstein \mf\/ $\,(M,g)$, acting 
on bivectors, has the same eigenvalues at every point\/ $\,x\in M$, then\/ 
$\,(M,g)\,$ is locally symmetric.
\endproclaim
This is immediate from the next result, proven in \S\rf\ (and, originally, in 
\cite{\hdg}):
\proclaim{Theorem \a\ri.2\ {\rm(\cite{\hdg}, p.~476, Theorem~7.1)}}If\/ 
$\,(M,g)\,$ is an oriented Riemannian Einstein four\mfd\/ and its self-dual 
curvature operator\/ $\,R^{\hs+}:\lap\to\lap\,$ has the same eigenvalues at 
every point, then\/ $\,R^{\hs+}$ is parallel.
\endproclaim

\head\S\fb. Further basics\endhead
Unless stated otherwise, all tensor fields are of class \ci\nh. For \of s 
$\,\xi,\eta$, \vf s $\,u,v,w\,$ and a \psr\ metric $\,g\,$ on any \mf, 
$\,\xi\wedge\eta,\hs d\hs\xi\,$ and the Lie derivative $\,\Cal L_wg\,$ are 
given by $\,(\xi\wedge\eta)(u,v)=\xi(u)\eta(v)-\xi(v)\eta(u)$, 
$\,(d\hs\xi)(u,v)=d_u[\xi(v)]-d_v[\xi(u)]-\xi([u,v])\,$ and, with 
$\,[\hskip2.5pt,\hskip1pt]\,$ denoting the Lie bracket,
$$(\Cal L_wg)(u,v)\,=\,d_w[g(u,v)]\,-\,g([u,w],v)\,-\,g(u,[v,w])\,.\ff\lw$$
On a \psr\ \mf\ $\,(M,g)\,$ we use the same symbol, such as $\,u\hs$, for a 
\vf\ and the corresponding \of\ $\,g(u,\,\cdot\,)$. Similarly, a vector-bundle 
morphism $\,\ah:\tm\to\tm\,$ is treated as a twice-contravariant tensor field, 
and as a twice-covariant one with $\,\ah(u,v)=g(\ah u,v)\,$ for \vf s $\,u,v$. 
In particular, a bi\vf\ $\,\ah\,$ (such as $\,v\wedge u$) is also regarded as 
a differential $\,2$-form, or a morphism $\,\ah:\tm\to\tm\,$ with 
$\,\ah^*=-\hs\ah\,$ (i.e., skew-adjoint at each point). Specifically, for 
bi\vf s $\,\ah,\ah'$ and \vf s $\,u,v,w$,
$$\aligned
\text{\rm a)}\hskip6pt&
v\wedge u\,=\,v\otimes u\,-\,u\otimes v\,,\hskip20pt
dw\,=\,\ef\,-\,\ef^*\nh,\hskip9pt\text{\rm where}\hskip7pt\ef=\nabla w\,,\\
\text{\rm b)}\hskip6pt&
(v\otimes u)w\,=\,\langle v,w\rangle u\,,\hskip20pt(v\wedge u)w\,=\,
\langle v,w\rangle u\,-\,\langle u,w\rangle v\,,\\
\text{\rm c)}\hskip6pt&
\langle\ah,v\wedge u\rangle\,=\,g(\ah v,\,u)\,,\hskip20pt
\langle\ah,\ah'\rangle\,=\,-\,\text{\rm Trace}\hskip3pt
(\ah\circ\ah')\hskip-.2pt/2\,,\\
\text{\rm d)}\hskip6pt&
\ah\circ(v\wedge u)=v\otimes(\ah u)-u\otimes(\ah v)\hs,\hskip4.5pt
\text{\rm e)}\hskip4.5pt\text{\rm Trace}\,[\ah\circ(v\wedge u)]
=-\hs2g(\ah v,u)\hs,\\
\text{\rm f)}\hskip6pt&
2\ef^*w\,=\,d\langle w,w\rangle\,,\qquad\text{\rm where}\quad
\ef\,=\,\nabla w\quad\text{\rm and}\quad\langle w,w\rangle=g(w,w)\,,
\endaligned\ff\vi$$
with $\,\langle v,w\rangle=g(v,w)$, $\,\langle u,w\rangle=g(u,w)\,$ in b). 
(Cf.\ \S\pr.) Here d) follows from b), and implies e), as 
$\,\hs\text{\rm Trace}\,(v\otimes u)=g(v,u)$, while $\,P=\nabla w:TM\to TM\,$ 
in a), f) acts by $\,\ef v=\nabla_{\!v}w$, and so 
$\,2\langle v,\ef^*w\rangle=2\langle\ef v,w\rangle=d_v\langle w,w\rangle$, 
which gives f).
\remark{Remark \a\fb.1}Let $\,(M,g)\,$ be a \psr\ Einstein \mf. Then 
$\,\hs\text{\rm div}\,\w=\,0\hs$. If, in addition, $\,M\,$ is oriented, 
$\,\dim M=4$, and the sign pattern of $\,g\,$ is $\,\pppp\,$ or $\,\mmpp$, 
then also $\,\hs\text{\rm div}\,\w^+=\,\hs\text{\rm div}\,\w^-=\,0\hs$.

In fact, these are well-known consequences of the second Bianchi identity 
(cf.\ \cite{\hdg}, pp.~460, 468). Here $\,\hs\text{\rm div}\,\ah$, for any 
covariant tensor field $\,\ah$, is the $\,g$-contraction of $\,\nabla\ah\,$ 
involving the first argument of $\,\ah\,$ and the differentiation argument.
\endremark
By a {\it \cvf\/} on a real \mf\ $\,M\,$ we mean a section $\,w\,$ of its 
complexified \tb. Sections of the ordinary (``real'') \tb\ of $\,M\,$ may be 
referred to as {\it real \vf s\/} on $\,M$. Thus, $\,w=u+iv\,$ with real \vf s 
$\,u=\,\text{\rm Re}\,w$, $\,v=\,\text{\rm Im}\,w$. Complex {\it bivector\/} 
fields are defined similarly.

All real-multilinear operations involving real vector/bi\vf s will, without 
further comment, be extended to complex vector (or, bivector) fields $\,v,w\,$ 
(or, $\,\ah,\ah'$), so as to become \cli\ in each argument. This includes the 
Lie bracket $\,[v,w]$, covariant derivative $\,\nabla_{\!v}w\,$ relative to 
any connection in the \tb, the inner product $\,g(v,w)$, the Lie derivative 
$\,\Cal L_wg\,$ for any given pseudo\-Riemannian metric $\,g$, the composite 
$\,\ah\circ\ah'$, as well as $\,\langle\ah,\ah'\rangle$, $\,\hs*\ah\,$ and 
$\,\ah v\,$ (cf.\ \f{\vi.c}). Note that $\,g(v,w),\hs\langle\ah,\ah'\rangle\,$ 
are {\it \cbi\/} (not sesquilinear!) in $\,v,w\,$ or $\,\ah,\ah'$, and a \ci\ 
\cvf\ $\,w\,$ is a \kf, i.e., $\,\Cal L_wg=0\hs$, if and only if its real and 
imaginary parts both are real \kf s.

Although the bi\vb\ $\,\biv$ of an oriented Lorentz\-i\-an four\mfd\ 
$\,(M,g)\,$ is a complex \vb\ with the multiplication by $\,\hs i\hs\,$ 
provided by $\,*\,$ (see \S\pr), it is also convenient to use the 
complexification $\,(\biv)^\bbC$ of its underlying real \vb. Then 
$\,(\biv)^\bbC=\lap\oplus\lam$, where $\,\lapm\,$ are, this time, the complex 
\vb s of fibre dimension $\,3$, obtained as the $\,(\pm\hs i)$-eigenspace 
bundles of $\,*\hs$. This is clear since $\,*\hs^2=-\,\text{\rm Id}\hs$, cf.\ 
\S\pr, and the complex-conjugation antiautomorphism 
$\,(\biv)^\bbC\to(\biv)^\bbC$ obviously sends $\,\lap\,$ onto $\,\lam$.

\head\S\ut. A unified treatment of all three cases\endhead
Throughout this section $\,(M,g)\,$ stands for a fixed oriented \psr\ 
four\mfd\ with a metric $\,g\,$ of one of the sign patterns \f{\sn}, while 
$\,\e$ is a complex vector bundle of of fibre dimension $\,3\,$ over 
$\,M$, and $\,\wy$ is a \cli\ bundle morphism $\,\e\to\e$. Our choices 
of $\,\e$ and $\,\wy$ are quite specific. Namely, when $\,g\,$ is 
Riemannian or neutral, $\,\e=[\lap]^\bbC$ is the complexification of the 
subbundle $\,\lap\,$ of $\,\biv$ (\S\pr) and $\,\wy$ is the unique 
$\,\bbC$-linear extension of $\,\w^+:\lap\to\lap\,$ to $\,[\lap]^\bbC$, while, 
if $\,g\,$ is Lo\-rentz\-i\-an, $\,\e=\lap\subset(\biv)^\bbC$ (see end of 
\S\fb) and $\,\wy$ is the restriction to $\,\e$ of the 
$\,\bbC$-linear extension of $\,\w:\biv\to\biv$ to $\,(\biv)^\bbC$. (The 
latter extension leaves $\,\e$ invariant, since $\,[\hs\w,\hs*\hs]=0\hs$, 
cf.\ \S\pr.)

We will use the symbol $\,\nabla\,$ for the connection in $\,\e\hs$ induced by 
the \lcc\ of $\,g$, and let $\,h\,$ stand for the \cbi\ fibre metric in 
$\,\e\hs$ which, in the Riemannian/neutral (or, Lo\-rentz\-i\-an) case is the 
unique \cbi\ extension of $\,\langle\,,\rangle\,$ (see \S\pr) from $\,\lap\,$ 
to $\,[\lap]^\bbC$ (or, respectively, the restriction to $\,\e=\lap\,$ of the 
\cbi\ extension of $\,\langle\,,\rangle\,$ from $\,\biv$ to $\,(\biv)^\bbC$). 
Note that, in all cases, $\,\nabla h=0\,$ and $\,\lap\,$ is a 
$\,\nabla$-parallel subbundle of $\,\biv$ or $\,(\biv)^\bbC$, since the \lcc\ 
of $\,g\,$ makes both $\,g\,$ and $\,*\,$ parallel. 

In this and the next two sections, the indices $\,j,k,l\,$ always vary in the 
range $\,\{1,2,3\}\,$ and repeated indices are summed over, unless explicitly 
stated otherwise.

Given $\,M,g,\e,\wy,\nabla,h\,$ as above, let us now fix any \ci\ local 
sections $\,\ah_j$ of $\,\e\hs$ which trivialize $\,\e\hs$ on an open set 
$\,\,U\subset M$. This gives rise to com\-plex-val\-ued functions $\,h_{jk}$ 
and \of s $\,\xi_j^k$ with
$$\alignedat2
&\text{\rm a)}\hskip4pt&&
\nabla\ah_j\nh=\xi_j^\ly\otimes\hskip1pt\ah_\ly\,,\hskip5pt\text{\rm i.e.,}
\hskip6pt\nabla_{\!v}\hskip1pt\ah_j\nh=\xi_j^\ly(v)\hskip1pt\ah_\ly\hskip6pt
\text{\rm for\ every\ tangent\ vector\ field}\hskip4ptv\hs,\\
&\text{\rm b)}\hskip4pt&&
dh_{jk}\nh=\xi_{jk}+\xi_{kj}\,,\hskip10pt
\text{\rm where}\hskip7pt\xi_{jk}\,=\,\xi_j^\ly h_{lk}\hskip5.5pt\text{\rm and}
\hskip6pth_{jk}\,=\,h(\ah_j,\hs\ah_k)\hs.\endalignedat\ff\nx$$
Thus, $\,h_{jk}$ are the component functions of the fibre metric $\,h\,$ and 
$\,\xi_j^k$ are the connection forms of $\,\nabla$, relative to the 
$\,\ah_j$, while \f{\nx.b} states that $\,\nabla h=0\hs$. For 
$\,\hs\text{\rm div}\hs\,$ as in Remark \a\fb.1 and all tangent vectors 
$\,v\,$ we have, with summation over $\,k$,
$$\aligned
&\text{\rm\ptmii i)}\quad[\nabla_{\!v}\w]\hskip1pt\ah_j\hs
=\,\theta_j^k(v)\hs\ah_k\,,\hskip40pt\text{\rm ii)}\quad
[\text{\rm div}\,\w]\hskip1pt\ah_j\,=\,\ah_k\,\theta_j^k\,,\hskip12pt\text{\rm 
where}\\
&\text{\rm iii)}\quad\w\nh\ah_j\,=\,\w_j^k\ah_k\hskip29pt\text{\rm and}
\hskip25pt\text{\rm iv)}\quad\theta_j^{\hs l}\,=\,d\hs\w_j^\ly\,
+\,\w_j^k\xi_k^{\hs l}\,-\,\w_k^\ly\xi_j^{\hs k}\hs.\endaligned\ff\nw$$
In fact, $\,\nabla_{\!v}$ applied to iii) gives i) (by \f{\nx.a}), and 
contracting i) we get ii) . (Here $\,\w\,$ might be replaced by $\,\wy$, as 
$\,\nabla\w^\pm$ are the $\,\lapm$ components of $\,\nabla\w$.)
\remark{Remark \a\ut.1}If $\,g\,$ is Riemannian or neutral, $\,\wy$ is 
the $\,\bbC$-linear extension of $\,\w^+$ to $\,[\lap]^\bbC$. Thus, in the 
Riemannian case, the eigenvalues of $\,\wy$ at every point are all real, as 
$\,\w:\biv\to\biv$ is self-adjoint.

If $\,g\,$ is Lo\-rentz\-i\-an, $\,\wy$ is, at each point, algebraically 
equivalent to $\,\w\,$ acting in $\,\biv$, since 
$\,\ah\mapsto\ah-i\hs[*\ah]\,$ is an isomorphism $\,\biv\to\lap\,$ of complex 
\vb s, sending $\,\w\,$ onto $\,\wy$ (as $\,[\hs\w,\hs*\hs]=0\hs$, cf.\ 
\S\pr).
\endremark

\head\S\cf. Calculations in a local orthonormal frame\endhead
As in \S\ut, the indices $\,j,k,l\,$ vary in the set $\,\{1,2,3\}$. The {\it 
Ricci symbol\/} $\,\ve_{jk\ly}$ will always stand for the signum of the 
permutation $\,(j,k,l)\,$ of $\,(1,2,3)$, if $\,j\ne k\ne l\ne j$, while 
$\,\ve_{jk\ly}=0\,$ if $\,j=k\,$ or $\,k=l\,$ or $\,l=j$. From now on we 
assume (cf.\ Remark \a\cf.3 below) that, for our $\,\ah_j$ and some 
$\,\ve_j\in\bbR\hs$,
$$\alignedat2
&\text{\rm\ptmi i)}\hskip8pt&&
\ve_j\ah_j\circ\ah_j=\,-\hs\text{\rm Id}\hskip8pt\text{\rm and}\hskip8pt
\ah_j\circ\ah_k\,=\,\ve_\ly\hskip1pt\ah_\ly\,=
\,-\,\ah_k\circ\ah_j\hskip9pt
\text{\rm if}\hskip6pt\ve_{jk\ly}=1\hs,\\
&\text{\rm ii)}\hskip8pt&&
\ve_1\ve_2\ve_3\,=\,1\hskip14pt\text{\rm and}\hskip14pt
\ve_j\hs\in\hs\{1,-1\}\,,\hskip5ptj=1,2,3\hs.\endalignedat\ff\ea$$
(No summing over $\,j,l$.) For a \cvf\ $\,w$, the \cvf s
$$v_j\,=\,\,\ah_jw\,,\qquad j=1,2,3\,,\ff\vw$$
satisfy, in view of \f{\ea} and skew-ad\-joint\-ness of the 
$\,\ah_j$, the relations
$$\alignedat2
&\text{\rm a)}\hskip5pt&&
\langle v_j,v_k\rangle=\hs\ve_j\langle w,w\rangle\,\delta_{jk}\hskip6pt
\text{\rm(no\ summation),}\hskip6pt\langle w,v_j\rangle=0\hs,\hskip6pt
j,k=1,2,3\hs,\\
&\text{\rm b)}\hskip5pt&&
\ah_jv_k=-\hs\ah_kv_j=\hs\ve_\ly v_\ly\hs,\hskip7pt
\ah_jv_j=-\hs\ve_jw\hskip8pt\text{\rm(no\ summing)}\hskip8pt
\text{if}\hskip7pt\ve_{jk\ly}=1\hs,\endalignedat\ff\av$$
with $\,\langle\,,\rangle\,$ standing for $\,g(\hskip3pt,\hskip2pt)$. From 
\f{\ea}, \f{\nx.b} and \f{\vi.c}, $\,h_{jk}=2\hskip.1pt\ve_j\hs\delta_{jk}$ 
(no summing), and so, again by \f{\nx.b}, the $\,\xi_{jk}$ are 
skew-sym\-met\-ric in $\,j,k$. Therefore, as $\,\ve_k\ve_\ly=\ve_j$ when 
$\,\{j,k,l\}=\{1,2,3\}\,$ (by \f{\ea.ii}), we have, from \f{\nx.a}, \f{\ea.ii},
$$\alignedat2
&\text{\rm\ptmi i)}\hskip8pt&&
\xi_j^{\hs j}\,=\,0\hskip4.2pt\text{\rm(no\ summing) \ and}\hskip5pt
\xi_j^{\hs k}\,=\,\ve_j\hskip1pt\xi_\ly\,,\hskip4.2pt
\xi_j^{\hs l}\,=\,-\hs\ve_j\hskip1pt\xi_k\hskip5pt
\text{\rm if}\hskip4.2pt\ve_{jk\ly}=1\,,\\
&\text{\rm ii)}\hskip8pt&&
\ve_j\nabla\ah_j\,=\,\xi_\ly\otimes\hskip1pt\ah_k\,
-\,\xi_k\otimes\hskip1pt\ah_\ly\qquad\text{\rm whenever}
\quad\ve_{jk\ly}=1\,,\endalignedat\ff\xj$$
with the \of s $\,\xi_j$ defined by $\,\xi_j=\ve_j\hskip1pt\xi_{kl}$ if 
$\,\ve_{jk\ly}=1$. Next, we define com\-plex-val\-ued functions 
$\,\lambda_j,\mu_j$, $\,j=1,2,3$, by
$$\lambda_j\,=\,\w_j^j\,,\quad\text{\rm and}\quad\mu_j\,=\,\ve_\ly\w_k^\ly
\quad\text{\rm if}\quad\{j,k,l\}\,=\,\{1,2,3\}\quad\text{\rm(no\ summing).}
\ff\lm$$
We always have $\,\hs\text{\rm Trace}\,\wy=0\,$ (\cite{\hdg}, p.~650); thus, 
for any function $\,\hs\text{\rm s}\hs\,$ on $\,\,U\nh$,
$$\alignedat2
&\text{\rm a)}\hskip5pt&&
\lambda_1\,+\,\lambda_2\,+\,\lambda_3\,=\,0\,,\\
&\text{\rm b)}\hskip5pt&&
L_1+L_2+L_3=0\hskip5pt\text{\rm if}\hskip4.5pt
L_j=(\lambda_k-\lambda_\ly)(\lambda_j+\,\text{\rm s}/12)\hskip4pt
\text{\rm whenever}\hskip4pt\ve_{jk\ly}=1\hs,\endalignedat\ff\ls$$
\remark{Remark \a\cf.1}By \f{\nw.iv}, \f{\xj.i} and \f{\lm}, 
$\,\theta_j^{\hs j}\nh=d\lambda_j\nh+2\mu_k\xi_k\nh-2\mu_\ly\xi_\ly$,
$\,\theta_j^{\hs k}\nh=\ve_k\hs d\mu_\ly\nh+\ve_j(\lambda_j\nh-\lambda_k)
\xi_\ly\nh+\ve_j\ve_k\mu_j\xi_k\nh-\mu_k\xi_j$, 
$\,\theta_j^{\hs l}\nh=\ve_\ly d\mu_k\nh+\ve_j(\lambda_\ly\nh-\lambda_j)\xi_k
\nh-\ve_j\ve_\ly\mu_j\xi_\ly\nh-\mu_\ly\xi_j$ (no summing), if 
$\,\ve_{jk\ly}\nh=1$. Hence, from \f{\nw.ii} and \f{\ea.i}, 
$\,[\text{\rm div}\,\w]\hskip1.5pt\ah_j=\hs\ah_j\left[\hs d\lambda_j+w_k
-w_\ly\right]$ whenever $\,\ve_{jk\ly}=1\,$ (no summing), for the \cvf s 
$\,w_j$ given by $\,w_j\,=\,2\hs\mu_j\xi_j\,+\,\ah_j\left[\hs d\mu_j\,+\,
\ve_k\ve_\ly(\lambda_k\,-\,\lambda_\ly)\xi_j\,+\,
\ve_k\mu_k\xi_\ly\,-\,\ve_\ly\mu_\ly\xi_k\right]$ if $\,\ve_{jk\ly}=1\,$ 
(no summing).

Consequently, if $\,\,\text{\rm div}\,\wy=0\,$ and the $\,\lambda_j$ are all 
constant, then there exists a \cvf\ $\,w\,$ with 
$\,w_1=w_2=w_3=w$. This is clear if one applies $\,\ah_j$ to the above formula 
for $\,[\text{\rm div}\,\w]\hskip1.5pt\ah_j$ and uses \f{\ea.i}.
\endremark
\remark{Remark \a\cf.2}For $\,M,g,\e,\wy,\hs U,\ah_j,\hs\ve_j,\xi_j$ as above, 
with \f{\ea},
\widestnumber\item{(iii)}\roster
\item"(i)"$dv_j\nh=\ve_j\hskip1pt\xi_\ly\wedge v_k\nh
-\ve_j\hskip1pt\xi_k\wedge v_\ly\nh+\ef_j$ if $\,\ve_{jk\ly}=1$, for 
$\,v_j,\ef_j$ given by \f{\vw} and 
$$\ef_j\,=\,\,\ah_j\circ\ef\,+\,\ef^*\nh\circ\,\ah_j\qquad\text{\rm for}\quad
\ef=\nabla w\,,\ff\pa$$
where $\,w\,$ is any given complex \ci\ \vf\ defined on $\,\,U$.
\item"(ii)"If $\,\w_j^k$ in \f{\nw.iii} are constant, the following two 
conditions are equivalent:
\itemitem"a)"$\nabla\wy=\,0\,$ everywhere in $\,\,U$.
\itemitem"b)"$\mu_1\xi_1=\mu_2\xi_2=\mu_3\xi_3$ and $\,(\lambda_k\,
-\,\lambda_\ly)\xi_j\,+\,\ve_\ly\mu_k\xi_\ly\,-\,\ve_k\mu_\ly\xi_k\,=\,0\,$ 
for $\,\lambda_j,\hs\mu_j$ given by \f{\lm} and any $\,j,k,l\,$ with 
$\,\{j,k,l\}\,=\,\{1,2,3\}$.
\item"(iii)"If $\,g\,$ is Einstein and $\,\ve_{jk\ly}=1$, then 
$\,\,d\hs\xi_j+\ve_j\hskip1pt\xi_k\wedge\xi_\ly\,=\,
-\hs(\w+\,\text{\rm s}/12)\hskip1pt\ah_j$.
\endroster
In fact, $\,\nabla v_j\nh=\ve_j\hskip1pt\xi_\ly\otimes v_k\nh
-\ve_j\hskip1pt\xi_k\otimes v_\ly+\,\ah_j\circ\ef\,$ if $\,\ve_{jk\ly}=1$, by 
\f{\xj.ii}, and so \f{\vi.a} yields (i).
Next, (ii) is clear from \f{\nw.i}, \f{\lm} and the formulae for 
$\,\theta_j^{\hs j},\theta_j^{\hs k},\theta_j^{\hs l}$ in Remark \a\cf.1.
Finally, let $\,g\,$ be Einstein. For fixed $\,j,k,\ly\,$ with 
$\,\ve_{jk\ly}=1$, \f{\cu.i}, \f{\xj.ii} and the formulae for 
$\,\xi\wedge\eta,\hs d\hs\xi\,$ in \S\fb\ give 
$\,-\hs\ve_k\ve_\ly\,\om/2=\hs d\hs\xi_j+\ve_j\hskip1pt\xi_k\wedge\xi_\ly$, 
where $\,\hs\om\hs\,$ is the com\-plex-val\-ued $\,2$-form with 
$\,\hs\om\hs(u,v)=h(R^\nabla(u,v)\ah_k,\ah_\ly)\,$ for any \vf s $\,u,v$, with 
$\,R^\nabla\,$ denoting the curvature of our connection 
$\,\nabla\,$ in $\,\e\hs$ (\S\ut). However, 
$\,R^\nabla(u,v)\ah_k=[R(u,v),\ah_k]$, where $\,[\hskip2.5pt,\hskip1pt]\,$ 
also stands for the commutator of bundle morphisms $\,\tmc\to\tmc$, and 
$\,R(u,v):\tmc\to\tmc$ is defined as in \f{\cu.i}; this is easily seen using 
\f{\cu.i} and the Leibniz-rule equality 
$\,\nabla_{\!u}\ah=[\nabla_{\!u},\ah]\,$ for such morphisms $\,\ah\,$ (with 
the commutator applied, this time, to operators acting on vector 
{\it fields}\hs). Hence, by Lemma~5.3 on p.~460 of \cite{\hdg}, 
$\,\hs\om\hs=R\hskip1pt[\ah_k,\ah_\ly]\,$ with 
$\,[\ah_k,\ah_\ly]=\ah_k\circ\ah_\ly-\ah_\ly\circ\ah_k$, i.e., as 
$\,[\ah_k,\hs\ah_\ly]=2\hskip.2pt\ve_j\hs\ah_j$ (by \f{\ea.i}) and 
$\,R=\w+\,\text{\rm s}/12\,$ (Remark \a\pr.1), we have 
$\,\hs\om\hs=2\hskip.2pt\ve_j\hs(\w+\,\text{\rm s}/12)\hs\ah_j$, and (iii) 
follows.
\endremark
\remark{Remark \a\cf.3}Let $\,M,g,\e,\wy$ be as in \S\ut. If  
$\,\wy(x):\e_x\to\e_x$ is \dia\ for every $\,x\in M\,$ and the set of its 
eigenvalues does not depend on $\,x$, then a suitable connected \nbd\ 
$\,\,U\,$ of any given point of $\,M\,$ admits \ci\ local trivializing 
sections $\,\ah_j$ of $\,\e$, $\,j=1,2,3$, satisfying conditions \f{\ea} along 
with \f{\xj.ii} for suitable $\,\ve_j,\xi_j$, and such that the corresponding 
com\-plex-val\-ued functions $\,\lambda_j,\mu_j$ in \f{\lm} are all constant, 
with $\,\mu_j=0$. Thus, $\,\w\nh\ah_j\,=\lambda_j\ah_j$ (no summing), for 
$\,j=1,2,3$, i.e., the $\,\lambda_j$ then are the (constant) eigenvalues of 
$\,\wy$.

Namely, by Lemma~6.15(ii),(iii) of \cite{\hdg}, p.~468, 
$\,\w\nh\ah_j\,=\lambda_j\ah_j$ for some \ci\ sections $\,\ah_j$ trivializing 
$\,\e\hs$ on such a set $\,\,U$, and constants $\,\lambda_j$. As $\,\wy$ is 
self-adjoint relative to $\,h$, cf.\ \S\pr, while $\,h\,$ is nondegenerate, 
the $\,\ah_j$ may be chosen so that $\,h_{jk}=2\hskip.1pt\ve_j\hs\delta_{jk}$ 
(no summing) with $\,\ve_j\in\{1,-\hs1\}$. Next, for sections $\,\ah,\beta\,$ 
of $\,\e$, the anticommutator $\,\{\ah,\beta\}=\ah\circ\beta+\beta\circ\ah\,$ 
equals $\,-\hs h(\ah,\beta)\,$ times $\,\hs\text{\rm Id}\hs$, and the 
commutator $\,[\ah,\beta]=\ah\circ\beta-\beta\circ\ah\,$ is a section of 
$\,\e$. (For $\,\{\ah,\beta\}\,$ one can verify this, in the Lo\-rentz\-i\-an 
case, using a basis of $\,\Lambda^{\!+\!}_xM$, $\,x\in M$, obtained by 
replacing each $\,\ah\,$ by $\,\ah-i\hs[*\ah]\,$ in a basis of $\,\biv$ of the 
form~(37.28) in \cite{\hdg}, p.~642; about the Riemannian and neutral cases, 
and for $\,[\ah,\beta]$, see \cite{\hdg}, formulae~(37.31), (37.29) on 
pp.~642, 643.) Thus, by \f{\vi.c}, 
$\,\ve_j\ah_j\circ\ah_j=\,-\hs\text{\rm Id}\hs$, $\,j=1,2,3$, and 
$\,\ah_j\circ\ah_k=\de_\ly\ve_\ly\hskip1pt\ah_\ly=-\,\ah_k\circ\ah_j$ 
whenever $\,\ve_{jk\ly}=1$, with some $\,\de_j\in\{1,-\hs1\}$, $\,j=1,2,3$. 
(In view of \f{\vi.c}, $\,\ah_j\circ\ah_k=[\ah_j,\ah_k]/2\,$ is 
$\,h$-or\-thog\-o\-nal to $\,\ah_j,\ah_k$.) Now, as 
$\,(\ah_j\circ\ah_k)\circ\ah_l=\ah_j\circ(\ah_k\circ\ah_l)$, we get 
$\,\de_\ly=\de_j$, and, similarly, $\,\de_1=\de_2=\de_3$. Applying an odd 
permutation to the $\,\ah_j$ and/or replacing them by $\,-\hs\ah_j$, if 
necessary, we now obtain \f{\ea}.
\endremark

\head\S\ms. The main structure theorem\endhead
The following result is a crucial step in our classification argument. We 
establish it using a refined version of the proof of Theorem 7.1 in 
\cite{\hdg} (pp. 477--479).
\proclaim{Theorem \a\ms.1}Suppose that\/ $\,(M,g)\,$ is an oriented \psr\ 
Einstein four\mfd\/ with one of the sign patterns\/ \f{\sn}, such that\/ 
$\,\wy:\e\to\e$, defined as in\/ {\rm\S\ut}, is \dia\ at every point 
and has constant eigenvalues.
\widestnumber\item{(ii)}\roster
\item"(i)"If\/ $\,g\,$ is positive definite, the self-dual Weyl tensor\/ 
$\,\w^+\,$ is parallel. 
\item"(ii)"If\/ $\,g\,$ is Lo\-rentz\-i\-an\/ $\,(\mppp)\,$ or neutral\/ 
$\,(\mmpp)\,$ and\/ $\,\wy$ is not parallel, then any given point of\/ 
$\,M\hs$ has a \nbd\/ $\,\,U\hs$ with\/ \ci\ \cvf s\/ $\,w,v_1,v_2,v_3$ 
which are linearly independent at every point of\/ $\,\,U\nh$, commute with 
every real \kf\ defined on any \os\ of\/ $\,\,U\nh$, and satisfy the 
inner-product and Lie-bracket relations
$$\aligned
&g(w,w)=g(v_j,v_j)=\gm\hskip5pt
\text{\rm(no summing),}\hskip5.5pt
g(w,v_j)=g(v_j,v_k)=0\hskip5pt
\text{\rm if}\hskip4.5ptj\ne k\hs,\\
&[w,v_j]=\rho_j\hs v_j\hskip10pt\text{\rm(no summing)}\,,\hskip12pt
[v_j,v_k]=0\hskip9pt\text{\rm for\ all}\hskip5pta,b\in\{1,2,3\}\hs,
\endaligned\ff\th$$
for some\/ $\,\gm\in\bbC\smallsetminus\{0\}$, where\/ $\,\rho_j\in\bbC\,$ are 
the three cubic roots of\/ $\,\gm^{\hs2}$, and both\/ 
$\,g,\hs[\hskip2.5pt,\hskip1pt]\,$ act \cbi ly on \cvf s.
\endroster
\endproclaim
\demo{Proof}Given $\,x\in M$, let us choose 
$\,\e,\wy,\hs U,\ah_j,\hs\ve_j,\xi_j,\lambda_j$ as in Remark \a\cf.3, with 
$\,x\in U$. Since $\,(M,g)\,$ is Einstein, $\,\,\text{\rm div}\,\wy=0\,$ 
(Remark \a\fb.1) and so, by Remark \a\cf.1, $\,w_1=w_2=w_3=w\,$ for some 
\cvf\ $\,w$, where the $\,w_j$ are as in Remark \a\cf.1 with 
$\,\mu_1=\mu_2=\mu_3=0\hs$. Now
\widestnumber\item{(b)}\roster
\item"(a)"$\,v_j=\hs(\lambda_\ly-\hs\lambda_k)\xi_j$ \ if 
$\,\hs\ve_{jk\ly}\nh=1$, for the \cvf s $\,\,v_j\nh=\ah_jw$,
\item"(b)"$\,d\hs\xi_j\,+\,\ve_j\hskip1pt\xi_k\wedge\xi_\ly\,=\,
-\hs(\lambda_j\,+\,\text{\rm s}/12)\hskip1pt\ah_j$ \ 
whenever\hskip7pt$\ve_{jk\ly}=1$,
\endroster
$\hs\text{\rm s}\hs\,$ being the scalar curvature; in fact, (a) follows if one 
applies $\,\ah_j$ to the formula for $\,w_j$ in Remark \a\cf.1 (with 
$\,\mu_j=0$ and $\,w_j=w$), using \f{\ea.i\r iii}, while (b) is obvious from 
Remark \a\cf.2(iii) with $\,\w\nh\ah_j\,=\lambda_j\ah_j$. We now define a 
constant $\,\phi\,$ by
$$\phi\,=\,(\lambda_j-\lambda_k)(\lambda_k-\lambda_\ly)(\lambda_\ly-\lambda_j)
\qquad\text{whenever}\quad\ve_{jk\ly}=1\,.\ff\fy$$
Throughout this proof we will write $\,\langle\,,\rangle\,$ instead of 
$\,g(\hskip3pt,\hskip2pt)$. For $\,\ef_j$ given by \f{\pa},
$$\alignedat2
&\text{\rm\ptmi i)}\hskip8pt&&
(\lambda_j-\lambda_k)(\lambda_j-\lambda_\ly)\hskip1pt\ef_j
\,+\,2\hskip1pt\ve_j(\lambda_k-\lambda_\ly)\,v_k\wedge v_\ly\,
=\,-\,(\lambda_j+\,\text{\rm s}/12)\,\phi\,\ah_j\hs,\hskip4pt\\
&\text{\rm ii)}\hskip8pt&&
(\lambda_j-\lambda_k)(\lambda_j-\lambda_\ly)\,\,\text{\rm div}\,w\,
+\,2\,(\lambda_k-\lambda_\ly)\langle w,w\rangle\,
=\,-\,(2\lambda_j+\,\text{\rm s}/6)\,\phi\,,\endalignedat\ff\dw$$
if $\,\ve_{jk\ly}=1$. In fact, since the $\,\lambda_j$ are constant, 
multiplying (b) above by $\,\phi\,$ and using (a) we obtain 
$\,(\lambda_j-\lambda_k)(\lambda_j-\lambda_\ly)\,dv_j+
\ve_j(\lambda_k-\lambda_\ly)\hs v_k\wedge v_\ly=
-\hs(\lambda_j+\,\text{\rm s}/12)\hs\phi\hs\ah_j$, if $\,\ve_{jk\ly}=1$. In 
view of Remark \a\cf.2(i), (a) and \f{\fy}, this is nothing else than 
\f{\dw.i}. Also, as 
$\,\hs\text{\rm div}\,w=\,\,\text{\rm Trace}_{\hskip.4pt\bbC}\nabla w$, taking 
the complex trace of the composites of both sides of \f{\dw.i} with $\,\ah_j$, 
we obtain \f{\dw.ii} from \f{\vi.e}, \f{\av} and \f{\ea.i}, since, by \f{\pa}, 
$\,2\,\text{\rm Trace}\,\ef=-\hs\ve_j\hs\text{\rm Trace}\,(\ah_j\circ\ef_j)\,$ 
(no summation). Next, for $\,\phi\,$ as in \f{\fy},
$$\phi\,=\,0\qquad\text{\rm if\ and\ only\ if}
\qquad\nabla\wy\,=\,0\quad\text{\rm identically.}\ff\fw$$
Namely, if $\,\phi=0\hs$, by \f{\fy}, (a) above, \f{\vw} and \f{\ea.i}, 
$\,w=0\,$ and $\,(\lambda_\ly-\lambda_k)\xi_j=0\,$ whenever 
$\,\ve_{jk\ly}=1$, and so (as $\,\mu_j=0\hs$, $\,j=1,2,3$), Remark \a\cf.2(ii) 
yields $\,\nabla\wy=0\hs$. Conversely, let $\,\nabla\wy=0\hs$. Remark 
\a\cf.2(ii) with $\,\mu_j=0\,$ now gives $\,(\lambda_\ly-\lambda_k)\xi_j=0\,$ 
whenever $\,\ve_{jk\ly}=1$. Hence $\,\phi=0\hs$, for if we had 
$\,\phi\ne0\hs$, the last relation and \f{\fy} would imply $\,\xi_j=0\hs$, 
$\,j=1,2,3$, i.e., from (b) above, 
$\,\lambda_1=\lambda_2=\lambda_3=-\,\text{\rm s}/12$, and, by \f{\fy}, 
$\,\phi\,$ would be zero anyway.

Since our assertion is immediate when $\,\nabla\wy=0\hs$, we now assume that
$$\phi\ne0\,,\qquad\text{\rm i.e.,}\quad
\lambda_1\ne\lambda_2\ne\lambda_3\ne\lambda_1\,.\ff\fn$$
(Cf.\ \f{\fy}, \f{\fw}.) We may treat \f{\dw.ii} as a system of three 
linear equations with two unknowns: $\,\,\text{\rm div}\,w\,$ and 
$\,\langle w,w\rangle$. This system's matrix has the $\,2\times2\,$ 
subdeterminants equal, by \f{\ls.a}, to 
$\,\pm\,6\lambda_j(\lambda_k-\lambda_\ly)^2$, $\,\ve_{jk\ly}=1$. They cannot 
be all zero, or else \f{\ls.a} would give 
$\,0=\lambda_j(\lambda_k-\lambda_\ly)=-\hs(\lambda_k+\lambda_\ly)
(\lambda_k-\lambda_\ly)=\lambda_\ly^2-\lambda_k^2$, if $\,\ve_{jk\ly}=1$, 
i.e., any two of the $\,\lambda_j$ would coincide up to a sign, so that, with 
the $\,\lambda_j$ suitably rearranged, $\,\lambda_2=\lambda_3=\pm\lambda_1$, 
contrary to \f{\fn}. The system \f{\dw.ii} thus has rank two, and can be 
solved for $\,\,\text{\rm div}\,w\,$ and $\,\langle w,w\rangle\,$ using 
determinants. Thus, $\,\langle w,w\rangle\,$ is constant since so are the 
coefficients of \f{\dw.ii} (cf.\ \f{\fy}); in addition, 
$\,\langle w,w\rangle\ne0$. Namely, if $\,\langle w,w\rangle\,$ were zero, we 
would have 
$\,\hs\text{\rm div}\,w=(\lambda_k-\lambda_\ly)(2\lambda_j+\,\text{\rm s}/6)$, 
$\,\ve_{jk\ly}=1$, by \f{\dw.ii}, \f{\fy} and \f{\fn}; summed over 
$\,j=1,2,3$, this would yield $\,\,\text{\rm div}\,w=0\,$ (cf.\ \f{\ls.b}); 
hence $\,(\lambda_k-\lambda_\ly)(2\lambda_j+\,\text{\rm s}/6)=0\hs$, 
$\,\ve_{jk\ly}=1$, which, in view of \f{\fn}, would imply that 
$\,2\lambda_j=-\,\text{\rm s}/6\,$ for $\,j=1,2,3$, contrary to \f{\fn}. 
Next,
$$\alignedat2
&\text{\rm\ptmi i)}\hskip8pt&&
\nabla_{\!v_j}w\,=\,\lambda_j(\lambda_k-\lambda_\ly)\,v_j\qquad\text{\rm 
whenever}\quad\ve_{jk\ly}=1\,,\\
&\text{\rm ii)}\hskip8pt&&
\text{\rm div}\,w\,=\,0\,,\hskip30pt\text{\rm iii)}\hskip8pt\text{\rm s}\,=\,0
\,,\hskip6pt\text{\rm i.e.,}\hskip5pt(M,g)\hskip5pt\text{\rm is\ Ricci-flat.}
\hskip30pt\endalignedat\ff\nd$$
In fact, both sides of \f{\dw.i} may be treated as bundle morphisms 
$\,\tmc\to\tmc$, and hence applied to the \cvf\ $\,v_j=\ah_jw$, giving, 
by \f{\pa}, \f{\av} and \f{\vi.b}, 
$\,(\lambda_j-\lambda_k)(\lambda_j-\lambda_\ly)\hskip1pt\ah_j\ef v_j
=\ve_j(\lambda_j+\,\text{\rm s}/12)\,\phi w$, whenever $\,\ve_{jk\ly}=1$. 
(Note that, from \f{\vi.b} and \f{\av.a}, $\,(v_k\wedge v_\ly)v_j=0\hs$, 
while, as $\,\langle w,w\rangle\,$ is constant, \f{\av.b} and \f{\vi.f} yield 
$\,\ef^*\nh\ah_jv_j=0\hs$.) Now, applying $\,\ah_j$ to both sides of the last 
equality, we obtain $\,(\lambda_j-\lambda_k)(\lambda_j-\lambda_\ly)\hskip1pt
\ef v_j=-\,(\lambda_j+\,\text{\rm s}/12)\phi v_j$ from \f{\ea.i} and \f{\vw}. 
Thus, since $\,\ef v=\nabla_{\!v}w\,$ for all vectors $\,v$, \f{\fy} and 
\f{\fn} imply that 
$\,\nabla_{\!v_j}w=(\lambda_k-\lambda_\ly)(\lambda_j+\,\text{\rm s}/12)\,v_j$ 
whenever $\,\ve_{jk\ly}=1$. Also, $\,w,v_1,v_2,v_3$ form an orthogonal 
trivialization of $\,[T\hs U]^\bbC$ (by \f{\av.a} with 
$\,\langle w,w\rangle\ne0$). Evaluating $\,\,\text{\rm div}\,w\,$ in that 
trivialization, we get, from \f{\av.a}, 
$\,\langle w,w\rangle\hskip2pt\text{\rm div}\,w=
\langle w,w\rangle\hskip2pt\text{\rm Trace}_{\hskip.4pt\bbC}\nabla w=
\sum_{j=1}^{\hs3}\ve_j\langle v_j\hs,\nabla_{\!v_j}w\rangle$, since 
$\,\langle w\hs,\nabla_{\!w}w\rangle=0\,$ as $\,\langle w,w\rangle\,$ is 
constant. Therefore, \f{\nd.ii} is immediate from the above formula for 
$\,\nabla_{\!v_j}w\,$ and \f{\av.a}, \f{\ls.b}. Finally, we have 
$\,(\lambda_k-\lambda_\ly)\langle w,w\rangle
=-\,(\lambda_j+\,\text{\rm s}/12)\,\phi$, $\,\ve_{jk\ly}=1$, from \f{\dw.ii} 
and \f{\nd.ii}. Summed over $\,j\,$ this gives $\,\,\text{\rm s}\hs\phi=0\,$ 
(by \f{\ls.a}), and so \f{\fn} yields \f{\nd.iii}, while \f{\nd.iii} 
and our formula for $\,\nabla_{\!v_j}w\,$ imply \f{\nd.i}.

Next, \f{\dw.ii} and \f{\nd.ii\r iii} give 
$\,\phi\lambda_j=(\lambda_\ly-\lambda_k)\langle w,w\rangle$, 
$\,\ve_{jk\ly}=1$, i.e., by \f{\fy} and \f{\fn}, $\,\langle w,w\rangle
=\lambda_j(\lambda_j-\lambda_k)(\lambda_j-\lambda_\ly)\,$ if 
$\,\ve_{jk\ly}=1$. 
However, this means that 
$\,\langle w,w\rangle=2\lambda_j^3+\lambda_1\lambda_2\lambda_3$, $\,j=1,2,3$. 
(Note that, whenever $\,\ve_{jk\ly}=1$, \f{\ls.a} yields 
$\,(\lambda_j-\lambda_k)(\lambda_j-\lambda_\ly)=
\lambda_j^2-\lambda_j(\lambda_k+\lambda_\ly)+\lambda_k\lambda_\ly=
2\lambda_j^2+\lambda_k\lambda_\ly$.) Thus, $\,\lambda_j^3=\mu\,$ for some 
complex number $\,\mu$, not depending on $\,j\in\{1,2,3\}\,$ and, by 
\f{\fn}, the $\,\lambda_j$ are the three cubic roots of $\,\mu$, so that 
$\,\lambda_1\lambda_2\lambda_3=\mu$. As 
$\,\langle w,w\rangle=2\lambda_j^3+\lambda_1\lambda_2\lambda_3$, we have 
$\,\lambda_j^3=-\hs\gm$, $\,j=1,2,3$, for the constant 
$\,\gm=-\hs\langle w,w\rangle/3\in\bbC\smallsetminus\{0\}$.

Hence, by \f{\fn}, the $\,\lambda_j$ cannot be all real. Thus, according to 
Remark \a\ut.1, \f{\fn} implies that $\,g\,$ cannot be Riemannian, i.e., 
$\,\phi=0\,$ in the Riemannmian case, which, in view of \f{\fw}, proves 
assertion (i).

Thus, we may assume that $\,g\,$ is Lo\-rentz\-i\-an or neutral and 
$\,\nabla\wy\ne0\hs$. Now
$$\aligned
&\text{\rm\ptmii i)}\hskip7pt\lambda_k\nh=z\lambda_j\hs,\hskip2.7pt
\lambda_\ly\,=\,\overline{z}\lambda_j\hskip6pt
\text{\rm if}\hskip4.5pt\ve_{jk\ly}=1\hs,\hskip7pt\text{\rm ii)}\hskip7pt
\lambda_k\,-\,\lambda_\ly\,=\,\pm\hs i\sqrt{3\,}\lambda_j\hs,\hskip6pt
\ve_{jk\ly}=1\hs,\\
&\text{\rm iii)}\hskip7pt[w,v_j]\,=\,\lambda_j(\lambda_\ly-\lambda_k)\,v_j
\hskip8pt\text{\rm if}\hskip6pt\ve_{jk\ly}=1\hs,\hskip11pt
\text{\rm iv)}\hskip7pt[v_j,v_k]\,=\,0\hskip7pt
\text{for\ all}\hskip6ptj,\,k\,,\endaligned\ff\lz$$
where $\,z=e^{\hs\pm2\pi i/3}$ for a suitable sign $\,\pm\,$ and 
$\,[\hskip2.5pt,\hskip1pt]\,$ is the Lie bracket. In fact, i) follows since 
$\,\lambda_j^3=-\hs\gm\ne0\hs$, while ii) is obvious from i) as 
$\,z-\overline z=\pm\hs i\sqrt3$. Next, applying \f{\dw.i} to $\,w\,$ and 
using \f{\pa}, \f{\vw}, \f{\vi.b}, \f{\av.a}, \f{\nd.ii}, 
\f{\fy} and \f{\fn}, we obtain the formula 
$\,\ah_j\nabla_{\!w}w=-\hs(\nabla w)^*v_j
+\lambda_j(\lambda_k-\lambda_\ly)\,v_j$, if $\,\ve_{jk\ly}=1$.
By \f{\nd.i}, $\,\langle(\nabla w)^*v_j,v_k\rangle=
\langle v_j,(\nabla w)v_k\rangle=0\,$ whenever $\,k\ne j\,$ (cf.\ 
\f{\av.a}). Since $\,w,v_1,v_2,v_3$ form a complex orthogonal basis at 
every point, our formula for $\,\ah_j\nabla_{\!w}w\,$ thus shows that, at each 
point, $\,\ah_j\nabla_{\!w}w\,$ is a combination of $\,w\,$ and $\,v_j$, i.e., 
by \f{\vw}, \f{\av.b} and \f{\ea.i}, $\,\nabla_{\!w}w=\psi_jw+\chi_jv_j$ 
for some functions $\,\psi_j$, $\,\chi_j$. As this is true for all 
$\,j\in\{1,2,3\}\,$ and $\,w\,$ does not depend on $\,j$, we have 
$\,\chi_j=0\hs$, while $\,\psi_j=0\,$ since $\,\langle w,w\rangle\,$ is 
constant. Consequently, $\,\nabla_{\!w}w=0\hs$. Furthermore, 
$\,\nabla_{\!w}\ah_j=0\,$ for all $\,j\,$ in view of \f{\xj.ii} and the 
relation $\,\langle\xi_j,w\rangle=0\,$ (immediate from (a) above, \f{\av.a} 
and \f{\fn}), and so \f{\vw} with $\,\nabla_{\!w}w=0\,$ gives 
$\,\nabla_{\!w}v_j=0\hs$, $\,j=1,2,3$. This, combined with \f{\nd.i} and the 
fact that $\,\nabla\,$ is torsionfree, proves \f{\lz.iii}. Next, since the 
$\,v_j$ are mutually orthogonal by \f{\av.a}, and every $\,\xi_j$ is a 
multiple of $\,v_j$ in view of (a) with \f{\fn}, we have, by  \f{\xj.ii}, 
$\,\nabla_{\!v_j}\ah_k=\ve_k\langle\xi_j,v_j\rangle\hs\ah_\ly\,$ (no 
summation) whenever $\,\ve_{jk\ly}=1$. Hence, by (a) and \f{\av.a}, 
$\,\nabla_{\!v_j}\ah_k
=\ve_j\ve_k(\lambda_\ly-\lambda_k)^{-1}\langle w,w\rangle\ah_\ly$, i.e., from 
\f{\ea.ii} and \f{\vw}, 
$\,[\nabla_{\!v_j}\ah_k]w
=\ve_\ly(\lambda_\ly-\lambda_k)^{-1}\langle w,w\rangle\hs v_\ly$, 
$\,\ve_{jk\ly}=1$. On the other hand, by \f{\nd.i} and \f{\av.b}, 
$\,\ah_k(\nabla_{\!v_j}w)
=-\hs\ve_\ly\lambda_j(\lambda_k-\lambda_\ly)\hs v_\ly$. From \f{\vw} and our 
expressions for $\,[\nabla_{\!v_j}\ah_k]w\,$ and 
$\,\ah_k(\nabla_{\!v_j}w)=0\,$ we now obtain 
$\,\nabla_{\!v_j}v_k=\nabla_{\!v_j}(\ah_kw)
=[\nabla_{\!v_j}\ah_k]w+\ah_k(\nabla_{\!v_j}w)=0\,$ if $\,\ve_{jk\ly}=1$, as 
\f{\lz.ii} with $\,\lambda_j^3=-\hs\gm=\langle w,w\rangle/3\,$ gives 
$\,(\lambda_\ly-\lambda_k)^{-1}\langle w,w\rangle
=\,\pm\,i\sqrt{3\,}\lambda_j^2$.) Similarly, $\,\nabla_{\!v_j}v_\ly=0\,$ 
if $\,\ve_{jk\ly}=1$. Thus, $\,\nabla_{\!v_j}v_k=0\,$ when $\,j\ne k$, proving 
\f{\lz.iv}.

Since, by \f{\av.a}, 
$\,\langle w,w\rangle=\langle v_j,v_j\rangle=-\,3\hs\gm\,$ (no summing) for 
$\,j=1,2,3$, the new \cvf s $\,\tilde w=\pm\,iw/\sqrt{3\,}\,$ and 
$\,\tilde v_j=iv_j/\sqrt{3\,}$, with the same sign $\,\pm\,$ as in 
\f{\lz.ii}, have $\,\langle\tilde w,\tilde w\rangle
=\langle\tilde v_j,\tilde v_j\rangle=\gm\,$ and are pairwise orthogonal by 
\f{\av.a}, while, from \f{\lz.ii} -- \f{\lz.iv}, 
$\,[\tilde v_j,\tilde v_k]=0\hs$, and 
$\,[\tilde w,\tilde v_j]=\lambda_j^2\tilde v_j$ (no summation). As 
$\,\lambda_j^3=-\hs\gm$, replacing $\,w,\,v_j$ with $\,\tilde w,\,\tilde v_j$ 
and setting $\,\rho_j=\lambda_j^2$, we now obtain (ii).

Finally, $\,w,v_j$ and $\,\tilde w,\tilde v_j$ commute with all \kf s 
since, up to permutations and sign changes, they are invariant under all 
isometries between connected \os s of $\,M$. Namely, by \f{\fn}, relations 
$\,\w\nh\ah_j\,=\lambda_j\ah_j$, \f{\ea.i}, \f{\xj.ii} and (a) above 
determine the $\,\ah_j,\xi_j,v_j$ and $\,w\,$ uniquely up to permutations and 
sign changes. This completes the proof.\quad\qed
\enddemo

\head\S\li. Complex \la s and real \mf s\endhead
Given a \rc\ \vs\ $\,\z\,$ of sections of a \rc\ vector bundle $\,\e\hs$ over 
a \mf\ $\,M$, we will say that $\,\z\,$ {\it trivializes} $\,\e\hs$ if it 
consists of \ci\ sections of $\,\e\hs$ and, for every $\,x\in M$, the 
evaluation operator $\,\psi\mapsto\psi(x)\,$ is an isomorphism $\,\z\to\e_x$. 
This amounts to requiring that $\,\dim\hs\z\,$ coincide with the fibre 
dimension of $\,\e\hs$ and each $\,v\in\z\,$ be either identically zero, or 
nonzero at every point of $\,M$. Equivalently, a basis of $\,\z\,$ then is a 
\ci\ trivialization of $\,\e$.

For instance, a simply transitive \la\ of \vf s\ on a \mf\ $\,M\,$ (see the 
appendix) is nothing else than a real \vs\ of \vf s on $\,M$, trivializing its 
(real) \tb, and closed under the Lie bracket.

Let the \rctb\ of a \mf\ $\,M\,$ be trivialized by a \rc\ \vs\ $\,\z\,$ of 
\rc\ \vf s on $\,M\,$ (cf.\ end of \S\fb). We will say that a \rc\ \vf\ 
$\,w\,$ defined of any \os\ $\,\,U\,$ of $\,M\,$ {\it commutes with\/} $\,\z$, 
and write $\,[w,\z\hs]=\{0\}$, if $\,w\,$ is of class \ci\ and $\,[w,v]=0\,$ 
for every\/ $\,v\in\z$. In view of the Jacobi identity, \rc\ \vf s $\,\z\,$ 
defined on a given open set $\,\,U\,$ and commuting with $\,\z\,$ form a \la.
\proclaim{Lemma \a\li.1}Let\/ $\,\z\,$ be a \rc\ \la\ of \rc\ \vf s on a real 
\mf\/ $\,M$, trivializing its \rctb. Then, the \rctb\ of any sufficiently 
small connected \nbd\/ $\,\,U\,$ of any given point\/ $\,x\,$ of\/ $\,M\,$ is 
trivialized by the \la\/ $\,\y\,$ of all\/ \rc\ \vf s defined on\/ $\,\,U\,$ 
and commuting with\/ $\,\z$.
\endproclaim
In fact, let $\,\hs\text{\rm D}\hs\,$ be the unique connection in the \rctb\ 
$\,\Cal T\,$ with $\,\hs\text{\rm D}_vw=[v,w]\,$ for all $\,v\in\z\,$ and all 
$\,C^1$ sections $\,w\,$ of $\,\Cal T$. Thus, $\,\hs\text{\rm D}\hs\,$ is 
flat: by \f{\cu.i}, $\,R^{\text{\rm D}}(v,w)u=[w,[v,u]]-[v,[w,u]]+[[v,w],u]\,$ 
whenever $\,v,w,u\in\z$, which is zero by the Jacobi identity. (As $\,\z\,$ is 
a \la, $\,[v,w]\in\z$, and so $\,\hs\text{\rm D}_{[v,w]}u=[[v,w],u]$.) Now 
$\,\y\,$ consists of all $\,\hs\text{\rm D}$-parallel sections of 
$\,\Cal T\hs$ on $\,\,U$,\quad\qed
\smallskip
For instance, the real \la\ $\,\x\,$ of left\inv\ \vf s on a Lie group $\,G\,$ 
trivializes its real \tb. A real \vf\ on an open connected subset $\,\,U\,$ of 
$\,G\,$ commutes with $\,\x\,$ if and only if it is the restriction to 
$\,\,U\,$ of a right\inv\ \vf\ on $\,G$. In fact, right\inv\ fields $\,w\,$ 
all commute with $\,\x$, since the flow of $\,w\,$ (or, of any $\,v\in\x$) 
consist of left (or, right) translations, while left and right translations 
commute due to associativity. The converse follows since both \la s are of 
dimension $\,\dim G\,$ (Lemma \a\li.1).
\remark{Remark \a\li.2}If a \rc\ \vs\ $\,\z\,$ of \rc\ \vf s on a \mf\ $\,M\,$ 
trivializes its \rc\ \tb, then any \rc\ \vf\ $\,w\,$ on $\,M\,$ with 
$\,[w,\z\hs]=\{0\}\,$ is a \rc\ \kf\ on $\,(M,g)\,$ for any \psr\ metric 
$\,g\,$ on $\,M\,$ such that $\,g(u,v)\,$ is constant whenever $\,u,v\in\z$. 
In fact, \f{\lw} then gives $\,(\Cal L_wg)(u,v)=0\,$ for all $\,u,v\in\z$.
\endremark
Let a complex \la\ $\,\z\,$ of \cvf s on a \mf\ $\,M\,$ trivialize its 
complexified \tb\ $\,\tmc$. We say that $\,\z\,$ {\it admits a real form\/} if 
$\,\hs\text{\rm Re}\,w\in\z\,$ for every $\,w\in\z$. This is obviously 
equivalent to the existence of a real \la\ $\,\x\,$ of real \vf s on $\,M$, 
trivializing its ordinary \tb\ $\,\tm$, and such that $\,\z=\x+i\x$, i.e., 
$\,\z\,$ is the complexification of $\,\x\,$ (or, $\,\x\,$ is a {\it real 
form\/} of $\,\z$). Clearly, $\,\x\,$ then is uniquely determined by $\,\z$, 
as $\,\x=\{\text{\rm Re}\,w:w\in\z\}=\{w\in\z:\,\text{\rm Im}\,w=0\}$. Thus, 
$\,\z\,$ admits a real form if and only if the real \vf s which are 
elements of $\,\z\,$ form a real \la\ trivializing $\,\tm$.
\remark{Remark \a\li.3}If a complex \la\ $\,\z\,$ of \cvf s on a \mf\ $\,M\,$ 
trivializes its complexified \tb\ and $\,\dimc\y=\dim M\,$ for the \la\ 
$\,\y\,$ of all \ci\ \cvf s $\,w\,$ on $\,M\,$ with $\,[w,\z\hs]=\{0\}$, then
\widestnumber\item{(ii)}\roster
\item"(i)"$\y\,$ trivializes the complexified \tb\ of $\,M$.
\item"(ii)"$\z\,$ admits a real form whenever $\,\y\,$ does.
\endroster
To see this, first note that Lemma \a\li.1 yields (i). Next, let $\,\tv\,$ be 
a real form of $\,\y$, and let a \cvf\ $\,w\,$ commute with $\,\y$, so that 
$\,[w,\y\hs]=\{0\}$. Since $\,\tv\subset\y$, we have $\,[w,\tv\hs]=\{0\}$. 
Therefore $\,[\hs\text{\rm Re}\,w,\tv\hs]=\{0\}$, as $\,\tv\,$ consists of 
real \vf s and $\,[\hskip2.5pt,\hskip1pt]\,$ is \cbi; this and 
relation $\,\y=\tv+i\tv\,$ now give $\,[\hs\text{\rm Re}\,w,\y\hs]=\{0\}$. The 
\la\ $\,\z'$ of all \cvf s commuting with $\,\y\,$ thus is closed under the 
real-part operator $\,\hs\text{\rm Re}\hs$. However, $\,\z\subset\z'$ and, by 
Lemma \a\li.1, $\,\dimc\z'\le\dim M=\dimc\z$, so that $\,\z'=\z$, which proves 
(ii).
\endremark
\proclaim{Lemma \a\li.4}Let\/ $\,\z\,$ be a complex \la\ of \cvf s on a 
\psrm\/ $\,(M,g)$, trivializing the complexified \tb\ of\/ $\,M\,$ and 
such that\/ $\,g(u,v)\,$ is constant for any\/ $\,u,v\in\z$. Then
\widestnumber\item{(b)}\roster
\item"(a)"$(M,g)\,$ is \lh.
\item"(b)"Under the additional assumption that\/ $\,[u,v]=0\,$ for every\/ 
$\,u\in\z\,$ and every real \kf\/ $\,v\,$ defined on any \os\ of\/ $\,M$, we 
have\/ $\,\hs\text{\rm Re}\,w\in\z\,$ whenever\/ $\,w\in\z$, i.e., $\,\z\,$ 
admits a real form.
\endroster
\endproclaim
In fact, by Lemma \a\li.1 and Remark \a\li.2, every vector in $\,T_xM$, 
$\,x\in M$, is the value at $\,x\,$ of some real \kf\ on a \nbd\ of $\,x$, 
which proves (a) (cf.\ \cite{\hdg}, p~546). Now let us fix $\,x\in M\,$ 
and choose $\,\,U,\y\,$ for $\,x,\z\,$ as in Lemma \a\li.1. Remark \a\li.2 
and our hypothesis show that $\,\y\,$ then is {\it precisely\/} the \la\ of 
all complex \kf s on $\,\,U$. Thus, $\,\y\,$ is closed under the real-part 
operator $\,\hs\text{\rm Re}\hs$, i.e., admits a real form, and Remark 
\a\li.3(ii) yields (b).
\quad\qed

\head\S\rf. Real forms of some specific complex \la s\endhead
We use the standard notation $\,\hs\text{\rm Ad}\hs\,$ for the adjoint 
representation of any given \la\ $\,\x$, so that 
$\,\hs\text{\rm Ad}\,v:\x\to\x\,$ is, for any $\,v\in\x$, given by 
$\,(\text{\rm Ad}\,v)w=[v,w]$.
\proclaim{Lemma \a\rf.1}Let a basis\/ $\,w,v_1,v_2,v_3$ of a four\diml\ 
complex \la\/ $\,\z\,$ satisfy conditions\/ \f{\th} for some \cbi\ symmetric 
form\/ $\,g\,$ on\/ $\,\z\,$ and a complex number\/ $\,\gm\ne0\hs$, where\/ 
$\,[\hskip2.5pt,\hskip1pt]\,$ is the Lie-algebra multiplication of\/ $\,\z\,$ 
and\/ $\,\rho_1,\rho_2,\rho_3$ are the three cubic roots of\/ $\,\gm^{\hs2}$. 
Also, let\/ $\,\x\subset\z\,$ be a four\diml\ real Lie subalgebra with\/ 
$\,\z=\x+i\x\,$ and\/ $\,g(\x,\x)\subset\bbR\hs$. In other words, $\,\x\,$ 
spans\/ $\,\z\,$ as a complex space and the form\/ $\,g\hs$ restricted to\/ 
$\,\x\hs$ is real-val\-ued.

Then\/ $\,w\in\x\,$ and there exist a three\diml\ real vector 
subspace\/ $\,V\hs$ of\/ $\,\x$, a linear operator\/ $\,\fe:V\to V$, and a 
real-val\-ued bilinear form\/ $\,\langle\,,\rangle\,$ on\/ $\,V$, satisfying 
conditions\/ \f{\bg} with\/ $\,u=|\gm\hs|^{-1/2}w\,$ and\/ 
$\,\da=\,\text{\rm sgn}\,\gm$, and such that\/ $\,\langle\,,\rangle,\fe\,$ 
are, for a suitable isomorphic identification\/ $\,V=\bbC\times\bbR\hs$, given 
by\/ {\rm(a),\hs(b)} in\/ {\rm\S\co} with some sign\/ $\,\pm\,$ and some\/ 
$\,\py\in\bbR\smallsetminus\{0\}$.
\endproclaim
\demo{Proof}We set $\,V=\x\cap\,\text{\rm Ker}\,\varPsi$ and define 
$\,\varPsi:\z\to\bbC\,$ to be the $\,\bbC$-linear functional with 
$\,\varPsi(w)=1\,$ and $\,\varPsi(v_j)=0\hs$, $\,j=1,2,3$. For any 
$\,u\in\z\smallsetminus\,\text{\rm Ker}\,\varPsi$,
\widestnumber\item{(b)}\roster
\item"(a)"$\hs\text{\rm Ad}\,u\,$ has the characteristic roots $\,0\,$ and 
$\,\varPsi(w)\hskip.8pt\rho_j$, $\,j=1,2,3$.
\item"(b)"$\dimr V=3\,$ and $\,\spanc V=\,\text{\rm Ker}\,\varPsi$. 
\endroster
In fact, by \f{\th}, $\,\,\text{\rm Ad}\,u:\z\to\z\,$ is \dia\, 
with the eigenvalues as in (a) for the eigenvectors $\,u\,$ and $\,v_j$, which 
proves (a). Also, as $\,\dimc\hs[\text{\rm Ker}\,\varPsi]=3$, 
our $\,\x\,$ cannot be contained in $\,\hs\text{\rm Ker}\,\varPsi$, and so the 
image $\,\varPsi(\x)\,$ is a nontrivial real vector subspace of $\,\bbC\hs$. 
For any fixed $\,u\in\x\smallsetminus\,\text{\rm Ker}\,\varPsi$, (a) gives 
$\,\varPsi(w)\hskip.8pt\rho_j\in\bbR\,$ for some $\,j\in\{1,2,3\}$. (In fact, 
as $\,\x\,$ contains a basis of $\,\z$, the characteristic roots of 
$\,\,\text{\rm Ad}\,w:\z\to\z\,$ coincide with those of 
$\,\,\text{\rm Ad}\,w:\x\to\x$, so that the number of nonreal ones among them 
is $\,0\,$ or $\,2$.) Thus, $\,\varPsi(\x)\,$ is contained 
in the union of the real lines $\,\bbR\hs\overline\rho_j\subset\bbC\hs$, 
$\,j=1,2,3$, i.e., must coincide with one of them, and we may fix 
$\,j\in\{1,2,3\}\,$ with $\,\varPsi(\x)\hs=\bbR\hs\overline{\rho_j}$. Now 
$\,\dimr V=3$, since $\,V=\x\cap\,\text{\rm Ker}\,\varPsi\,$ is the kernel of 
$\,\varPsi:\x\to\bbC\hs$. Also, $\,\x\,$ spans $\,\z$, so that vectors in 
$\,\,\x$, linearly independent over $\,\bbR\hs$, are also linearly independent 
over $\,\bbC\,$ in $\,\z$. This implies (b): 
$\,\hs\spanc V=\,\text{\rm Ker}\,\varPsi\,$ as 
$\,\hs\spanc V\subset\,\text{\rm Ker}\,\varPsi\,$ and 
$\,\dimc\hs[\hs\spanc V]=\dimc\hs[\hs\text{\rm Ker}\,\varPsi]=3$.

As $\,\dimr V=3$, we may choose $\,u\in\x\smallsetminus\{0\}\,$ which is 
$\,g$-or\-thog\-o\-nal to $\,V$. By (b), $\,u\,$ then is also 
$\,g$-or\-thog\-o\-nal to $\,\,\text{\rm Ker}\,\varPsi$. Hence, in view of 
\f{\th}, $\,u\in\bbC w$, i.e., $\,u=\varPsi(u)w\,$ with $\,\varPsi(u)\ne0\hs$. 
Also, $\,\varPsi(u)\hskip.8pt\rho_j$ is real, for $\,j\,$ chosen above (as 
$\,\varPsi(u)\in\bbR\hs\overline{\rho_j}\hs$), and hence so is its cube 
$\,[\varPsi(u)]^3\gm^{\hs2}$. On the other hand, \f{\th} gives 
$\,[\varPsi(u)]^2\gm=g(u,u)\in\bbR\hs$. Consequently, the numbers 
$\,\varPsi(u)\hs\gm$, $\,\varPsi(u)$, $\,\rho_j$ and $\,\gm\,$ are all real, 
while $\,w\in\x$, as $\,\x\,$ contains $\,u=\varPsi(u)w\,$ and 
$\,\varPsi(u)\in\bbR\smallsetminus\{0\}$.

Since $\,\gm\in\bbR\smallsetminus\{0\}$, replacing such $\,u\,$ by 
$\,|\gm\hs|^{-1/2}w\,$ and letting $\,\langle\,,\rangle\,$ stand for the 
restriction of $\,g\,$ to $\,V$, we now obtain $\,\langle u,u\rangle=\da\,$ 
with $\,\da=\,\text{\rm sgn}\,\gm\in\{1,-\hs1\}$. 

The real $\,3$-space $\,V=\x\cap\,\text{\rm Ker}\,\varPsi\,$ is 
$\,(\text{\rm Ad}\,u)\hs$\inv, since so are $\,\x\,$ (as $\,u\in\x$) and 
$\,\hs\text{\rm Ker}\,\varPsi\,$ (by \f{\th} with $\,u=|\gm\hs|^{-1/2}w$). 
The restriction $\,\fe:V\to V\hs$ of $\,\hs\text{\rm Ad}\,u\,$ is 
self-adjoint, since that is the case for $\,\fe,V,\langle\,,\rangle\,$ 
replaced by $\,\hs\text{\rm Ad}\,u,\z,g\,$ (as 
$\,\hs\text{\rm Ad}\,u:\z\to\z\,$ is diagonalized by the 
$\,g$-or\-thog\-o\-nal basis $\,w,v_1,v_2,v_3$, cf.\ \f{\th}). Combining 
(a) with our assumptions about the cubes $\,\rho_j^3$ and the fact that 
$\,\dimr V=3\,$ is odd, we see that $\,\fe\,$ has the characteristic roots 
$\,\py,\hs\py\qe,\hs\py\overline \qe$, where $\,\qe=e^{2\pi i/3}\,$ and 
$\,\py\in\bbR\smallsetminus\{0\}$, and we may choose 
$\,\xi,\eta,\zeta\in V\hs$ such that $\,\zeta\,$ and $\,\xi+i\eta\,$ are 
eigenvectors of $\,\hs\text{\rm Ad}\,u:\z\to\z\,$ for the eigenvalues 
$\,\py\,$ and $\,\py\qe$. (Since $\,\py\qe\notin\bbR\hs$, this implies that 
$\,\xi,\eta\,$ are linearly independent over $\,\bbR$.) By \f{\th}, 
$\,\zeta\,$ and $\,\xi+i\eta\,$ are complex multiples of $\,v_j,v_k$ for some 
$\,j,k$. Thus, $\,g(\zeta,\zeta)\ne0\hs$, i.e., $\,\zeta\,$ may be normalized 
so that $\,\langle\zeta,\zeta\rangle=\pm\hs1\,$ for some sign $\,\pm\hs$, 
while $\,g(\zeta,\hs\xi+i\eta)=0\hs$, and so 
$\,\langle\zeta,\xi\rangle=\langle\zeta,\eta\rangle=0\hs$, as $\,g\,$ is 
real-val\-ued on $\,V\subset\x$. Next, 
$\,\langle \fe\xi,\eta\rangle=\langle\xi,\fe\eta\rangle\,$ since $\,\fe\,$ is 
self-adjoint, so that $\,\langle\xi,\xi\rangle+\langle\eta,\eta\rangle=0\,$ in 
view of the eigenvector relation $\,\fe\xi+i\fe\eta=\py\qe(\xi+i\eta)\,$ with 
$\,\qe=(\sqrt{3\,}i-1)/2$. Finally, let $\,c\,$ be a complex number with 
$\,2\hs\overline c^{\hs2}=-\hs g(\xi+i\eta,\xi+i\eta)$. Thus, $\,c\ne0\hs$, 
since $\,g(v_k,v_k)\ne0\hs$, and it is easy to verify that the isomorphism 
$\,V\to\bbC\times\bbR\,$ sending the basis $\,\xi,\eta,\zeta\,$ onto 
$\,(c,0),(-\hs ic,0),(0,1)\,$ has the required properties. This completes the 
proof.\quad\qed
\enddemo
\demo{Proofs of Theorems \a\cl.1, \a\cn.1 and \a\ri.2}In all three cases, 
$\,R-\w\,$ is a constant multiple of the identity (Remark \a\pr.1), and so the 
hypotheses of Theorem \a\ms.1 are satisfied (cf.\ Remark \a\ut.1). If $\,g\,$ 
is Riemannian, Theorem \a\ms.1(i) yields  Theorem \a\ri.2. If $\,g\,$ is 
Lo\-rentz\-i\-an and $\,\nabla\w=0\hs$, i.e., $\,\nabla R=0\hs$, Theorem 41.5 
of \cite{\hdg} (pp.~662--663) implies (a) or (b) in Theorem \a\cl.1, as the 
\diy\ condition excludes option (c) in \cite{\hdg} on p.~663. The only 
remaining cases now are those named in (ii) of Theorem \a\ms.1, the conclusion 
of which shows that Lemma \a\rf.1 can be applied to the \la\ 
$\,\z=\hs\spanc\hs\{w,v_1,v_2,v_3\}\,$ and its real form $\,\x\,$ which exists 
in view of Lemma \a\li.4(b). As a result, $\,(M,g)\,$ is obtained as in 
Example \a\co.2(i) or (ii); the situation where $\,\da=-\hs1\,$ and $\,\pm\,$ 
is $\,-\,$ cannot occur, as it would lead to the sign pattern $\,\mmmp\hs$, 
which is not one of \f{\sn}.\quad\qed
\enddemo

\head Appendix. Simply transitive \la s of \vf s\endhead
In this section we prove Corollary \a\tr.3 which, although well-known, seems 
to lack a convenient reference; we need it for a conclusion in Example \a\co.2.

A {\it simply transitive \la\ of \vf s\/} on a \mf\ $\,M\,$ is any \vs\ 
$\,\x\,$ of \ci\ (real) \vf s on $\,M$, closed under the Lie bracket and such 
that the evaluation operator $\,\x\ni w\mapsto w(x)\in T_xM\,$ is bijective 
for every $\,x\in M$. An example is the \la\ of left\inv\ \vf s on a Lie group.

Given a simply transitive \la\ $\,\x\,$ of \vf s on a \mf\ $\,M\,$ and a fixed 
point $\,y\in M$, the {\it exponential mapping\/} $\,\xp:U_y\to M\,$ {\it 
for\/} $\,\x$, {\it centered at\/} $\,y$, is given by $\,\xp(v)=x(1)$, where 
$\,\,U_y$ is the set of all $\,v\in\x\,$ for which an integral curve 
$\,t\mapsto x(t)\,$ of $\,v\,$ with $\,x(0)=y\,$ can be defined on the whole 
interval $\,[0,1]$. It is clear that $\,\,U_y$ is a \nbd\ of $\,0\,$ in 
$\,\x\,$ and, for every $\,v\in U_y$ and $\,t\in[\hs0,1\hs]$, we have 
$\,tv\in U_y$ and $\,x(t)=\xp(tv)$, with $\,x(t)\,$ as above.

Let $\,Q:\bbC\to\bbC\,$ be the entire function with $\,Q(z)=(1-e^{-z})/z\,$ if 
$\,z\ne0\,$ and $\,Q(0)=1$. Its Maclaurin series defines $\,Q(A)\,$ for any 
linear operator $\,A:V\to V\,$ in a \vs\ $\,V\hs$ with $\,\dim V<\infty$. 
Thus, with $\,\hs\text{\rm Ad}\hs\,$ as in \S\rf, $\,Q(\text{\rm Ad}\,v)\,
=\,\sum_{\hs k=0}^{\hs\infty}\hs(-\hs\text{\rm Ad}\,v)^k/[\hs(k+1)!\hs]\,$ for 
a \la\ $\,\x\,$ with $\,\dim\x<\infty\,$ and $\,v\in\x$.
\proclaim{Proposition \a\tr.1}Let\/ $\,\x\,$ be a simply transitive \la\ of 
\vf s on a \mf\/ $\,M$, and let\/ $\,d\xp_v:\x\to T_{\xp(v)}M\,$ be the 
differential at\/ $\,v\in U_y$ of the exponential mapping of\/ $\,\x\,$ 
centered at a point\/ $\,y\in M$, with the usual identification\/ 
$\,T_v\x=\x$. Then\/ $\,d\xp_v$ equals the composite mapping in which\/ 
$\,Q(\text{\rm Ad}\,v):\x\to\x$, defined above, is followed by the evaluation 
isomorphism\/ $\,\x\to T_{\xp(v)}M$.
\endproclaim
\demo{Proof}For any $\,C^\infty$ mapping $\,(s,t)\mapsto x(s,t)\in M\,$ 
of a rectangle $\,K\subset\rto$, let $\,u_s,u_t:K\to\x\,$ assign to 
$\,(s,t)\,$ the unique elements of $\,\x\,$ which coincide, at $\,x(s,t)$, 
with $\,\partial x/\partial s\,$ and, respectively, 
$\,\partial x/\partial t\,$ (that is, with the velocity at $\,s$, or $\,t$, of 
the curve $\,s\mapsto x(s,t)\,$ or $\,t\mapsto x(s,t)$). Using subscripts for 
partial derivatives of $\,u_s,u_t$ we thus have 
$\,u_{st},u_{ts},u_{stt}:K\to\x\,$ with $\,u_{st}=\partial u_s/\partial t$, 
etc.; we also let $\,[u_s,u_t]:K\to\x\,$ stand for the valuewise bracket of 
the Lie-algebra valued functions $\,u_s,u_t$. In local coordinates $\,x^j$ at 
any given $\,x_0=x(s_0,t_0)$, the vector fields $\,u_s(s,t),u_t(s,t)\,$ have 
some component functions $\,u_s^j(s,t,x),u_t^j(s,t,x)$, also depending on a 
point $\,x\,$ near $\,x_0$. Thus, 
$\,u_s^j(s,t,x(s,t))=\hs\partial\hs[x^j(s,t)]/\partial s\,$ and 
$\,u_t^j(s,t,x(s,t))=\hs\partial\hs[x^j(s,t)]/\partial t$. Applying 
$\,\partial/\partial t\,$ to the first relation, $\,\partial/\partial s\,$ to 
the second, and using equality of mixed partial derivatives for the 
$\,x^j(s,t)$, we get $\,\partial u_s^j/\partial t-\hs\partial u_t^j/\partial s
=u_s^k\hs\partial_ku_t^j-u_t^k\hs\partial_ku_s^j$, with 
$\,\partial_k=\hs\partial/\partial x^k$, which is the coordinate form of the 
identity $\,u_{st}-u_{ts}=[u_s,u_t]$. If $\,u_{tt}=0\,$ for all 
$\,(s,t)\in K$, taking $\,\partial/\partial t\,$ of that identity, we obtain 
the {\it Jacobi equation\/} $\,u_{stt}=[u_{st},\hs u_t]\,$ (as 
$\,u_{tst}=u_{tts}=0$).

It is clear that $\,u_{tt}=0\,$ identically if and only if 
$\,t\mapsto x(s,t)\,$ is, for each fixed $\,s$, an integral curve of some \vf\ 
$\,v(s)\in\x$. Then, obviously, $\,u_t(s,t)=v(s)$. 

Now let $\,u_{tt}=0\,$ for all $\,(s,t)$, and let $\,K\,$ intersect the 
$\,s$-axis $\,\bbR\times\{0\}$. The Jacobi equation (see above) reads 
$\,\partial u_{st}/\partial t=-\hs[\text{Ad}\,v(s)]\,u_{st}$, with 
$\,v(s)=u_t(s,t)$, and so $\,u_{st}(s,t)=e^{-\hs t\,\text{Ad}\,v(s)}\,w(s)$, 
where $\,w(s)=u_{st}(s,0)$. Since 
$\,d\hs[t\hs Q(t\,\text{\rm Ad}\,v)]/dt=e^{-\hs t\,\text{Ad}\,v}\,$ (cf.\ our 
formula for\hskip3pt$Q(\text{\rm Ad}\,v)$), we get 
$\,u_s(s,t)=u_s(s,0)+t\hs Q(t\,\text{\rm Ad}\,v(s))\hs w(s)$, as both sides 
satisfy the same initial value problem in the variable $\,t$.

Finally, let $\,K=I\times[\hs0,1\hs]\,$ and $\,x(s,t)=\xp(tv(s))\,$ for some 
interval $\,I\,$ and some $\,C^\infty$ curve $\,I\ni s\mapsto v(s)\in U_y$. 
Thus, $\,u_{tt}=0\,$ identically and $\,u_t(s,t)=v(s)$, so that 
$\,u_{ts}(s,t)=\dot v(s)$, with $\,\dot v=\hs dv/ds$. Also, $\,x(s,0)=y$, and 
hence $\,u_s(s,0)=0\hs$. Evaluating at $\,(s,0)\,$ the identity 
$\,u_{st}-u_{ts}=[u_s,u_t]$, established above, and setting 
$\,w(s)=u_{st}(s,0)\,$ as in the preceding paragraph, we thus get 
$\,w(s)=u_{ts}(s,0)=\dot v(s)$. Writing $\,v,\dot v\,$ instead od 
$\,v(s),\hs dv/ds\,$ we now see that $\,u_s(s,1)\,$ equals the pre\-im\-age of 
$\,d\xp_v\dot v\,$ under the evaluation isomorphism $\,\x\to T_{\xp(v)}M\,$ 
(cf.\ the definition of $\,u_s$) while 
$\,u_s(s,1)=Q(\text{\rm Ad}\,v)\hs\dot v$, as one sees setting $\,t=1\,$ in 
$\,u_s(s,t)=u_s(s,0)+t\hs Q(t\,\text{\rm Ad}\,v(s))\hs w(s)$. This completes 
the proof.\quad\qed
\enddemo
\proclaim{Corollary \a\tr.2}Given a simply transitive \la\/ $\,\x\,$ of \vf s 
on a \mf\/ $\,M\,$ and a point\/ $\,y\in M\,$ there exists a neighborhhod\/ 
$\,\,U\,$ of\/ $\,0\,$ in\/ $\,\x$ such that\/ $\,\,U\subset\,U_y$ and the 
exponential mapping\/ $\,\xp:U_y\to M\,$ sends\/ $\,\,U\,$ diffeomorphically 
onto an open subset of\/ $\,M$. For any\/ $\,\,U\,$ with this property, 
$\,Q(\text{\rm Ad}\,v):\x\to\x\,$ is an isomorphism for every\/ $\,v\in U$, 
and the pull\-back under\/ $\,\xp\,$ of any \vf\/ $\,w\in\x\,$ is the \vf\ 
on\/ $\,\,U\,$ given by\/ $\,\,U\ni v\,\mapsto\,[Q(\text{\rm Ad}\,v)]^{-1}w$.
\endproclaim
In fact, $\,Q(\text{\rm Ad}\,v)\,$ is an isomorphism by Proposition \a\tr.1, 
since $\,d\xp_v$ is.\quad\qed
\smallskip
By Corollary \a\tr.2, the local diffeomorphism type of a simply transitive 
\la\ of \vf s is determined by its Lie-algebra isomorphism type. Since every 
\fdi\ \la\ is the \la\ of some Lie group, this yields
\proclaim{Corollary \a\tr.3}Given a simply transitive \la\/ $\,\x\,$ of \vf s 
on a \mf\/ $\,M$, there exists a Lie group\/ $\,G\,$ with the following 
property\/{\rm:} Every point of\/ $\,M\,$ has a \nbd\/ $\,\,U\,$ which 
may be diffeomorphically identified with an open set\/ $\,\,U'\subset G\,$ so 
as to make\/ $\,\x\,$ restricted to\/ $\,\,U\,$ appear as the \la\ of the 
restrictions to\/ $\,\,U'$ of all left\inv\ \vf s on\/ $\,G$.\quad\qed
\endproclaim



\enddocument